\documentclass[10pt]{article}
\usepackage{amsmath, amsthm, amssymb}
\usepackage[mathscr]{eucal}

\newcommand{\A}{\mathbb{A}}
\newcommand{\B}{\mathbb{B}}
\newcommand{\C}{\mathbb{C}}

\newcommand{\G}{\mathbb{G}}

\newcommand{\N}{\mathbb{N}}

\newcommand{\Q}{\mathbb{Q}}
\newcommand{\R}{\mathbb{R}}

\newcommand{\T}{\mathbb{T}}
\newcommand{\Z}{\mathbb{Z}}

\newcommand{\cA}{\mathcal{A}}

\newcommand{\cH}{\mathcal{H}}

\newcommand{\cV}{\mathcal{V}}

\newcommand{\scA}{\mathscr{A}}
\newcommand{\scE}{\mathscr{E}}
\newcommand{\scM}{\mathscr{M}}

\def\gg{\mathfrak{g}}

\newcommand{\gp}{\mathfrak{p}}
\newcommand{\gq}{\mathfrak{q}}
\newcommand{\gt}{\mathfrak{t}}

\newcommand{\Ab}{\mathbf{Ab}}
\newcommand{\cl}{\mathrm{class.}}
\newcommand{\coker}{\mathrm{coker}\ }
\newcommand{\Cong}{\mathrm{Cong}}

\newcommand{\cts}{\mathrm{cts}}
\newcommand{\disc}{\mathrm{disc.}}
\newcommand{\Ext}{\mathrm{Ext}}

\newcommand{\Group}{\mathrm{Group}}
\newcommand{\GL}{\mathrm{GL}}
\newcommand{\Hom}{\mathrm{Hom}}
\newcommand{\ind}{\mathrm{ind}}

\newcommand{\la}{\mathrm{loc.an.}}
\newcommand{\lalg}{\mathrm{loc.alg.}}
\newcommand{\Jacq}{\mathrm{Jacq}}
\newcommand{\Lie}{\mathrm{Lie}}
\newcommand{\limprojs}{\lim_{{\leftarrow}\atop{s}}\ }
\newcommand{\limprojsd}{\left(\lim_{\leftarrow\atop s}\right)^{1}}
\newcommand{\limdKp}{\lim_{\to\atop K_{\gp}}\ }

\newcommand{\rank}{\mathrm{rank}}
\newcommand{\ram}{\mathrm{ramified}}
\newcommand{\relLie}{\mathrm{rel.Lie}}

\newcommand{\SL}{\mathrm{SL}}
\newcommand{\Spec}{\mathrm{Spec}\ }
\newcommand{\sll}{\mathfrak{sl}}
\newcommand{\smooth}{\mathrm{smooth}}
\newcommand{\SO}{\mathrm{SO}}
\newcommand{\sph}{\mathrm{sph}}
\newcommand{\Spin}{\mathrm{Spin}}
\newcommand{\Sp}{\mathrm{Sp}}
\newcommand{\sq}{\hfill$\square$}
\newcommand{\st}{\mathrm{st}}

\newcommand{\SU}{\mathrm{SU}}

\newcommand{\unram}{\mathrm{unramified}}

\newtheorem{theorem}{Theorem}

\newtheorem{lemma}{Lemma}
\newtheorem{corollary}{Corollary}
\newtheorem{proposition}{Proposition}

\newtheorem*{thm}{Theorem}

\theoremstyle{definition}
\newtheorem{definition}{Definition}

\theoremstyle{remark}
\newtheorem*{remark}{Remark}
\newtheorem*{acknowledgement}{Acknowledgements}

\newtheorem*{question}{Question}

\begin{document}

\title{Construction of eigenvarieties
in small cohomological dimensions
for semi-simple, simply connected groups}
\author{Richard Hill}
\maketitle

\tableofcontents

\begin{abstract}
We study low order terms of Emerton's spectral sequence
for simply connected, simple groups.
As a result, for real rank 1 groups, we show that Emerton's method
for constructing eigenvarieties is successful
in cohomological dimension 1.
For real rank 2 groups, we show that a slight modification of
Emerton's method allows one to construct
eigenvarieties in cohomological dimension 2.
\end{abstract}

Throughout this paper we shall use the following standard notation:
\begin{itemize}
	\item
	$k$ is an algebraic number field, fixed throughout.
	\item
	$\gp$, $\gq$ denote finite primes of $k$,
	and $k_{\gp}$, $k_{\gq}$ the corresponding local fields.
	\item
	$k_{\infty}=k\otimes_{\Q}\R$ is the product of the
	archimedean completions of $k$.
	\item
	$\A$ is the ad\`ele ring of $k$.
	\item
	$\A_{f}$ is the ring of finite ad\`eles of $k$.
	\item
	For a finite set $S$ of places of $k$, we let
	$$
		k_{S}
		=
		\prod_{v\in S} k_{v},
		\qquad
		\A^{S}
		=
		{\prod_{v\notin S}}' k_{v}.
	$$
\end{itemize}

\section{Introduction and Statements of Results}

\subsection{Interpolation of classical automorphic representations}

Let $\G$ be a connected, algebraically simply connected, semi-simple
 group over a number field $k$.
We fix once and for all
 a maximal compact subgroup $K_{\infty}\subset\G(k_{\infty})$.
Our assumptions on $\G$ imply that $K_{\infty}$ is connected
in the archimedean topology.
This paper is concerned with the cohomology of the
following ``Shimura manifolds'':
$$
	Y(K_{f})
	=
	\G(k) \backslash \G(\A) / K_{\infty}K_{f},
$$
where $K_{f}$ is a compact open subgroup of $\G(\A_{f})$.
Let $W$ be an irreducible finite dimensional algebraic representation of
 $\G$ over a field extension $E/k$.
Such a representation gives rise to a local system
 $\cV_{W}$ on $Y(K_{f})$.
We shall refer to the cohomology groups of this local system
as the ``classical cohomology groups'':
$$
	H^{\bullet}_{\cl}(K_{f},W)
	:=
	H^{\bullet}(Y(K_{f}),\cV_{W}).
$$
It is convenient to consider the direct limit over all levels $K_{f}$
 of these cohomology groups:
$$
	H^{\bullet}_{\cl}(\G,W)
	=
	\lim_{\to\atop K_{f}}
	H^{\bullet}_{\cl}(K_{f},W).
$$
There is a smooth action
 of $\G(\A_{f})$ on $H^{\bullet}_{\cl}(\G,W)$, and we may
 recover the finite level cohomology groups as spaces of $K_{f}$-invariants:
$$
	H^{\bullet}_{\cl}(K_{f},W)
	=
	H^{\bullet}_{\cl}(\G,W)^{K_{f}}.
$$
It has become clear that only a very restricted
 class of smooth representations of $\G(\A_{f})$
 may occur as subquotients of the classical cohomology $H^{n}_{\cl}(\G,W)$.
For example, in the case $E=\C$,
 Ramanujan's Conjecture (Deligne's Theorem) gives an
 archimedian bound on the eigenvalues of the Hecke operators.
We shall be concerned with non-archimedean bounds.

Fix once and for all a finite prime $\gp$ of $k$
 over which $\G$ is quasi-split.
Fix a Borel subgroup $\B$ of $\G\times_{k}k_{\gp}$
and a maximal torus $\T\subset \B$.
We let $E$ be a finite extension of $k_{\gp}$, large enough
so that $\G$ splits over $E$.
It follows that the irreducible algebraic representations of $\G$ over $E$
 are absolutely irreducible (\S24.5 of \cite{borelbook}).
By the highest weight theorem (\S24.3 of \cite{borelbook}),
 an irreducible representation $W$ of $\G$ over $E$
 is determined by its highest weight $\psi_{W}$,
 which is an algebraic character
 $\psi_{W}:\T\times_{k_{\gp}} E\to \GL_{1}/E$.

By a \emph{tame level} we shall mean a compact open subgroup
$K^{\gp}\subset \G(\A_{f}^{\gp})$.
Fix a tame level $K^{\gp}$, and consider the spaces of $K^{\gp}$-invariants:
$$
	H^{\bullet}_{\cl}(K^{\gp},W)
	=
	H^{\bullet}_{\cl}(\G,W)^{K^{\gp}}.
$$
The group $\G(k_{\gp})$ acts smoothly on $H^{\bullet}_{\cl}(K^{\gp},W)$.
We also have commuting actions of the level $K^{\gp}$ Hecke algebra:
$$
	\cH(K^{\gp})
	:=
	\left\{
	f : K^{\gp}\backslash \G(\A_{f}^{\gp})/K^{\gp}\to E :
	\hbox{$f$ has compact support}
	\right\}.
$$
In order to describe the representations of $\cH(K^{\gp})$,
recall the tensor product decomposition:
\begin{equation}
	\label{tensor}
 	\cH(K^{\gp})
	=
	\cH(K^{\gp})^{\ram}
	\otimes
	\cH(K^{\gp})^{\sph},
\end{equation}
where $\cH(K^{\gp})^{\sph}$ is commutative but infinitely generated, and
$\cH(K^{\gp})^{\ram}$ is non-commutative but finitely generated.
Consequently the irreducible representations of $\cH(K^{\gp})$ are finite-dimensional.

Let $\gq\ne \gp$ be a finite prime of $k$.
We shall say that $\gq$ is unramified in $K^{\gp}$ if
\begin{itemize}
	\item[(a)]
	$\G$ is quasi-split over $k_{\gq}$, and splits over an
	unramified extension of $k_{\gq}$, and
	\item[(b)]
	$K^{\gp}\cap \G(k_{\gq})$ is a hyper-special
	maximal compact subgroup of $\G(k_{\gq})$
	(see \cite{tits}).
\end{itemize}
Let $S$ be the set of finite primes $\gq\ne \gp$, which are ramified in $K^{\gp}$.
This is a finite set, and we have
$$
	K^{\gp}
	=
	K_{S}\times \prod_{\gq \; \unram} K_{\gq},
	\qquad
	K_{S}=K^{\gp}\cap \G(k_{S}),
	\quad
	K_{\gq}
	=
	K^{\gp}\cap\G(k_{\gq}).
$$
This gives the tensor product decomposition (\ref{tensor}),
 where we take
$$
	\cH(K^{\gp})^{\ram}
	=
	\cH(K_{S}),\qquad
	\cH(K^{\gp})^{\sph}
	=
	{\bigotimes_{\gq \; \unram}}'
	\cH(K_{\gp}).
$$
For each unramified prime $\gq$, the Satake isomorphism
 (Theorem 4.1 of \cite{cartier}) shows that $\cH(K_{\gq})$
 is finitely generated and commutative.
Hence the irreducible representations of $\cH(K^{\gp})^{\sph}$
 over $\bar E$ are $1$-dimensional,
 and may be identified with elements of $(\Spec \cH(K^{\gp})^{\sph})(\bar E)$.
Since the global Hecke algebra is infinitely generated,
 $\Spec \cH(K^{\gp})^{\sph}$ is an infinite dimensional space.
One might expect that the representations which occur as subquotients of
$H^{\bullet}_{\cl}(K^{\gp},W)$ are evenly spread around this space.
There is an increasing body of evidence that this is not the case,
and that in fact these representations are contained in a finite dimensional subset of
 $\Spec \cH(K^{\gp})^{\sph}$, independent of $W$.

More precisely, let $\pi$ be an irreducible representation of
$\G(k_{\gp})\times \cH(K^{\gp})$, which occurs as a subquotient of
$H^{n}_{\cl}(K^{\gp}, W)\otimes_{E}\bar E$.
We may decompose $\pi$ as a tensor product:
$$
	\pi
	=
	\pi_{\gp}
	\otimes
	\pi^{\ram}
	\otimes
	\pi^{\sph},
$$
where $\pi^{\sph}$ is a character of $\cH(K^{\gp})^{\sph}$;
$\pi^{\ram}$ is an irreducible representation of $\cH(K^{\gp})^{\ram}$ and
$\pi_{\gp}$ is an irreducible smooth representation of $\G(k_{\gp})$.
We can say very little about the pair $(W,\pi)$ in this generality,
so we shall make another restriction.
We shall write $\Jacq_{\B}(\pi_{\gp})$ for the Jacquet module of $\pi_{\gp}$,
with respect to $\B(k_{\gp})$.
The Jacquet module is a smooth, finite dimensional representation of $\T(k_{\gp})$.
It seems possible to say something about those pairs $(\pi,W)$
 for which $\pi_{\gp}$ has non-zero Jacquet module.
Such representations $\pi_{\gp}$ are also said to have \emph{finite slope}.
Classically for $\GL_{2}/\Q$, representations of finite slope
 correspond to Hecke eigenforms for which the eigenvalue of $U_{\gp}$ is non-zero.
By Frobenius reciprocity, such a $\pi_{\gp}$ is a submodule of a smoothly induced
 representation $\ind_{\B(k_{\gp})}^{\G(k_{\gp})}\theta$,
 where $\theta:\T(k_{\gp})\to \bar E^{\times}$ is a smooth character.
In order to combine the highest weight $\psi_{W}$, which is an algebraic character
of $\T$, and the smooth character $\theta$ of $\T(k_{\gp})$,
we introduce the following rigid analytic space (see \cite{schneider} for background
in rigid analytic geometry):
$$
	\hat T(A)
	=
	\Hom_{k_{\gp}-\la}(\T(k_{\gp}),A^{\times}),
	\quad
	A \hbox{ a commutative Banach algebra over $E$}.
$$
Emerton defined the \emph{classical point}
 corresponding to $\pi$ to be the pair
$$
	(\theta\psi_{W},\pi^{\sph})
	\in
	\left(\hat T \times \Spec\cH(K^{\gp})^{\sph}\right)(\bar E).
$$
We let $E(n,K^{\gp})_{\cl}$ denote the set of all classical points.
Emerton defined the \emph{eigenvariety} $E(n,K^{\gp})$
 to be the rigid analytic Zariski closure of $E(n,K^{\gp})_{\cl}$
 in $\hat T\times \Spec\cH(K^{\gp})^{\sph}$.

Concretely, this means that for every unramified prime
$\gq$ and each generator $T^{i}_{\gq}$ for the
Hecke algebra $\cH(K_{\gq})$, there is a holomorphic
function $t_{\gq}^{i}$ on $E(n,K^{\gp})$ such that for every
representation $\pi$ in $H^{n}_{\cl}(K^{\gp},W)\otimes \bar E$
of finite slope at $\gp$,
the action of $T_{\gq}^{i}$ on $\pi$ is by scalar multiplication
by $t_{\gq}^{i}(x)$, where $x$ is the corresponding classical point.

One also obtains a description of the action of the ramified part
 of the Hecke algebra.
This description is different, since irreducible
 representations of $\cH(K^{\gp})^{\ram}$ are
 finite dimensional rather than $1$-dimensional.
Instead one finds that there is a coherent sheaf $\scM$
 of $\cH(K^{\gp})^{\ram}$-modules over $E(n , K^{\gp})$,
 such that, roughly speaking, the action of $\cH(K^{\gp})^{\ram}$ on the fibre of
 a classical point describes the action of $\cH(K^{\gp})^{\ram}$ on the
 corresponding part of the classical cohomology.
A precise statement is given in Theorem \ref{emertonmain} below.

Emerton introduced a criterion
 (Definition \ref{emertoncrit} below),
 according to which the Eigencurve $E(n,K^{\gp})$
 is finite dimensional.
More precisely, he was able to prove that
the projection $E(n,K^{\gp})\to \hat T$ is finite.
If we let $\gt$ denote the Lie algebra of $\T(\bar E)$, then there is
a map given by differentiation at the identity element:
$$
	\hat T
	\to
	\check \gt,
$$
where $\check \gt$ is the dual space of $\gt$.
It is worth noting that the image in $\check \gt$ of a classical point
depends only on the highest weight $\psi_{W}$, since smooth
characters have zero derivative.
Emerton also proved, assuming his criterion,
 that the projection $E(n,K^{\gp})\to \check\gt$ has discrete fibres.
As a result, one knows that the dimension of the eigencurve is at
most the absolute rank of $\G$.

The purpose of this paper is to investigate Emerton's criterion for
connected, simply connected, simple groups.
Specifically, we show that Emerton's criterion holds for all such groups
in dimension $n=1$.
Emerton's criterion typically fails in dimension $n=2$.
However we prove a weaker form of the criterion for $n=2$,
and we show that the weaker criterion is sufficient for most purposes.

\subsection{Emerton's Criterion}

Let $p$ be the rational prime below $\gp$.
In \cite{emerton} Emerton introduced the following $p$-adic Banach spaces:
$$
	\tilde H^{\bullet}(K^{\gp},\Q_{p})
	=
	\left(
	\limprojs
	\lim_{\to\atop K_{\gp}}
	H^{\bullet}(Y(K_{\gp}K^{\gp}),\Z/p^{s})
	\right)
	\otimes_{\Z_{p}} \Q_{p}.
$$
For convenience, we also consider the direct limits of these spaces over
all tame levels $K^{\gp}$:
$$
	\tilde H^{\bullet}(\G,\Q_{p})
	=
	\lim_{\to \atop K^{\gp}}
	\tilde H^{\bullet}(K^{\gp},\Q_{p}).
$$
We have the following actions on these spaces:
\begin{itemize}
	\item
	The group $\G(\A_{f}^{\gp})$ acts smoothly on $\tilde H^{\bullet}(\G,\Q_{p})$;
	the subspace $\tilde H^{\bullet}(K^{\gp},\Q_{p})$ may be recovered as the
	$K^{\gp}$-invariants:
	$$
		\tilde H^{\bullet}(K^{\gp},\Q_{p})
		=
		\tilde H^{\bullet}(\G,\Q_{p})^{K^{\gp}}.
	$$
	\item
	The Hecke algebra $\cH(K^{\gp})$
	acts on $\tilde H^{\bullet}(K^{\gp},\Q_{p})\otimes E$.
	\item
	The group $\G(k_{\gp})$ acts continuously,
	but not usually smoothly on the Banach space
	$\tilde H^{\bullet}(K^{\gp},\Q_{\gp})$.
	This is an admissible continuous representation of $\G(k_{\gp})$
	in the sense of \cite{schneider-teitelbaum}
	(or \cite{emertonAnalytic}, Definition 7.2.1).
	\item
	Recall that we have fixed a finite extension
	$E/k_{\gp}$, over which $\G$ splits.
	We let
	$$
		\tilde H^{\bullet}(K^{\gp},E)
		=
		\tilde H^{\bullet}(K^{\gp},\Q_{p})\otimes_{\Q_{p}}E.
	$$
	The group $\G(k_{\gp})$ is a $\gp$-adic analytic group.
	Hence, we may define the subspace of $k_{\gp}$-locally analytic vectors
	in $\tilde H^{\bullet}(K^{\gp},E)$ (see \cite{emertonAnalytic}):
	$$
		\tilde H^{\bullet}(K^{\gp},E)_{\la}.
	$$
	This subspace is $\G(k_{\gp})$-invariant,
	 and is an admissible locally analytic representation of $\G(k_{\gp})$
	 (in the sense of \cite{emertonAnalytic}, Definition 7.2.7).
	The Lie algebra $\gg$ of $\G$
	also acts on the subspace $\tilde H^{\bullet}(K^{\gp},E)_{\la}$.
\end{itemize}

For an irreducible algebraic representation $W$ of $\G$ over $E$,
we shall write $\check W$ be the contragredient representation.
Emerton showed (Theorem 2.2.11 of \cite{emerton})
 that there is a spectral sequence:
\begin{equation}
	\label{specseq}
	E_{2}^{p,q}
	=
	\Ext^{p}_{\gg}(\check W, \tilde H^{q}(K^{\gp},E)_{\la})
	\implies
	H^{p+q}_{\cl}(K^{\gp},W).
\end{equation}
Taking the direct limit over the tame levels $K^{\gp}$, there is also a
spectral sequence
(Theorem 0.5 of \cite{emerton}):
\begin{equation}
	\label{specseqlimit}
	\Ext^{p}_{\gg}(\check W, \tilde H^{q}(\G,E)_{\la})
	\implies
	H^{p+q}_{\cl}(\G,W).
\end{equation}
In particular, there is an edge map
\begin{equation}
	\label{edge}
	H^{n}_{\cl}(\G,W)
	\to
	\Hom_{\gg}(\check W, \tilde H^{n}(\G,E)_{\la}).
\end{equation}

\begin{definition}
	\label{emertoncrit}
	We shall say that $\G$ satisfies \emph{Emerton's criterion in dimension $n$}
	if the following holds:
	\begin{center}
		For every $W$, the edge map (\ref{edge}) is an isomorphism.
	\end{center}
	This is equivalent to the edge maps
	$H^{n}_{\cl}(K^{\gp},W)\to\Hom_{\gg}(\check W, \tilde H^{n}(K^{\gp},E)_{\la})$
	being isomorphisms for every $W$ and every tame level $K^{\gp}$.
\end{definition}

\begin{theorem}[Theorem 0.7 of \cite{emerton}]
	\label{emertonmain}
	Suppose Emerton's criterion holds for $\G$ in dimension $n$.
	Then we have:
	\begin{enumerate}
		\item
		Projection onto the first factor induces a finite map
		$E(n,K^{\gp})\to\hat{T}$.
		\item
		The map $E(n,K^{\gp})\to\check{\gt}$ has discrete fibres.
		\item
		If $(\chi,\lambda)$ is a point of the Eigencurve such that $\chi$ is
		locally algebraic and of non-critical slope
		(in the sense of \cite{emertonJacq}, Definition 4.4.3),
		then $(\chi,\lambda)$ is a classical point.
		\item
		There is a coherent sheaf $\scM$ of $\cH(K^{\gp})^{\ram}$-modules
		 over $E(n , K^{\gp})$ with the following property.
		For any classical point $(\theta\psi_{W},\lambda)\in E(n,K^{\gp})$
		 of non-critical slope,
		 the fibre of $\scM$ over the point $(\theta\psi_{W}, \lambda)$
		 is isomorphic (as a $\cH(K^{\gp})^{\ram}$-module) to 
		 the dual of the $(\theta\psi_{W},\lambda)$-eigenspace
		 of the Jacquet module of $H^{n}_{\cl}(K^{\gp}, \check W )$.
	\end{enumerate}
\end{theorem}

In fact Emerton proved this theorem for all reductive groups $\G/k$.
He verified his criterion in the case $\G=\GL_{2}/\Q$, $n=1$.
He also pointed out that the criterion always holds for $n=0$,
 since the edge map at $(0,0)$ for any first quadrant $E_{2}^{\bullet,\bullet}$
 spectral sequence is an isomorphism.
Of course the cohomology of $\G$ is usually uninteresting
in dimension $0$, but his argument can be applied in the case
where the derived subgroup of $\G$ has real rank zero.
This is the case, for example, when $\G$ is a torus, or the multiplicative group of a definite quaternion algebra.

\subsection{Our Main Results}

For our main results, $\G$ is connected, simple and algebraically
simply connected.
We shall also assume that $\G(k_{\infty})$ is not compact.
We do not need to assume that $\G$ is absolutely simple.
We shall prove the following.

\begin{theorem}
	\label{crit1}
	Emerton's criterion holds in dimension $1$.
\end{theorem}

For cohomological dimensions $2$ and higher,
Emerton's criterion is quite rare.
We shall instead use the following criterion.

\begin{definition}
	We shall say that $\G$ satisfies the \emph{weak Emerton criterion}
	in dimension $n$ if
	\begin{itemize}
		\item[(a)]
		for every non-trivial irreducible $W$,
		the edge map (\ref{edge}) is an isomorphism, and
		\item[(b)]
		for the trivial representation $W$,
		the edge map (\ref{edge}) is injective,
		and its cokernel is a finite dimensional trivial representation of $\G(\A_{f})$.
	\end{itemize}
\end{definition}

By simple modifications to Emerton's proof of
Theorem \ref{emertonmain}, we shall prove the following in \S\ref{variationproof}.

\begin{theorem}
	\label{variation}
	If the weak Emerton criterion holds for $\G$ in dimension $n$,
	then
	\begin{enumerate}
		\item
		Projection onto the first factor induces a finite map
		$E(n,K^{\gp})\to\hat{T}$.
		\item
		The map $E(n,K^{\gp})\to\check{\gt}$ has discrete fibres.
		\item
		If $(\chi,\lambda)$ is a point of the Eigencurve such that $\chi$ is
		locally algebraic and of non-critical slope,
		then either $(\chi,\lambda)$ is a classical point
		or $(\chi,\lambda)$ is the trivial representation
		of $\T(k_{\gp})\times\cH(K^{\gp})^{\sph}$.
	\end{enumerate}
\end{theorem}

In order to state our next theorems,
we recall the definition of the congruence kernel.
As before, $\G/k$ is simple, connected and simply connected
 and $\G(k_{\infty})$ is not compact.
By a \emph{congruence subgroup} of $\G(k)$, we shall mean a subgroup of the form
$$
	\Gamma(K_{f})
	=
	\G(k)
	\cap
	(\G(k_{\infty})\times K_{f}),
$$
where $K_{f}\subset\G(\A_{f})$ is compact and open.
Any two congruence subgroups are commensurable.

An \emph{arithmetic subgroup} is a subgroup of $\G(k)$,
 which is commensurable with a congruence subgroup.
In particular, every congruence subgroup is arithmetic.
The \emph{congruence subgroup problem}
(see the survey articles \cite{raghunathan,rapinchuk})
 is the problem of determining
 the difference between arithmetic subgroups and congruence subgroups.
In particular, one could naively ask whether every arithmetic subgroup of $\G$
is a congruence subgroup.
In order to study this question more precisely,
Serre introduced two completions of $\G(k)$:
$$
	\hat{\G}(k)
	=
	\lim_{\leftarrow \atop K_{f}}
	\G(k) / \Gamma(K_{f}),
$$
$$
	\tilde{\G}(k)
	=
	\lim_{\leftarrow \atop \Gamma\; \mathrm{arithmetic}}
	\G(k) / \Gamma.
$$
There is a continuous surjective group homomorphism
$\tilde{\G}(k) \to \hat{\G}(k)$.
The \emph{congruence kernel} $\Cong(\G)$ is defined to be the
kernel of this map.
Recall the following:

\begin{theorem}[Strong Approximation Theorem
	\cite{kneser1,kneser-hasseprinciple,kneser,platonov,prasad}]
	\label{SAT}
	Suppose $\G/k$ is connected, simple, and algebraically simply connected.
	Let $S$ be a set of places of $k$, such that $\G(k_{S})$ is not compact.
	Then $\G(k)\G(k_{S})$ is dense in $\G(\A)$.
\end{theorem}

Under our assumptions on $\G$, the strong approximation theorem
implies that $\hat{\G}(k)=\G(\A_{f})$,
and we have the following extension of topological groups:
$$
	1 \to \Cong(\G) \to \tilde{\G}(k) \to \G(\A_{f}) \to 1.
$$
By the \emph{real rank} of $\G$, we shall mean the sum
$$
	m
	=
	\sum_{\nu |\infty}\rank_{k_{\nu}}\G.
$$
It follows from the non-compactness of $\G(k_{\infty})$,
that the real rank of $\G$ is at least $1$.
Serre \cite{serre3} has conjectured that for $\G$ simple, simply connected and of
real rank at least $2$, the congruence kernel is finite;
 for real rank $1$ groups he conjectured that the congruence kernel is infinite.
These conjectures have been proved in many cases
and there are no proven counterexamples
(see the surveys \cite{raghunathan,rapinchuk}).

Our next result is the following.

\begin{theorem}
	\label{crit2}
	If the congruence kernel of $\G$ is finite then
	the weak Emerton criterion holds in dimension $2$.
\end{theorem}

Theorems \ref{crit1} and \ref{crit2}
follow from our main auxiliary results:

\begin{theorem}
	\label{mainA}
	Let $\G$ be as described above.
	Then $\tilde H^{0}(\G,E)=E$, with the trivial action of $\G(\A_{f})$.
\end{theorem}

\begin{theorem}
	\label{mainB}
	Let $\G$ be as described above.
	Then
	$$
		\tilde H^{1}(\G,E)
		=
		\Hom_{\cts}(\Cong(\G),E)_{\G(\A_{f}^{\gp})-\smooth},
	$$
	where $\Cong(\G)$ denotes the congruence kernel of $\G$.
\end{theorem}

The reduction of Theorem \ref{crit1} to Theorem \ref{mainA}
is given in \S\ref{crit1proof}, and the reduction of Theorem \ref{crit2}
to Theorem \ref{mainB} is given in \S\ref{crit2proof}.
Theorem \ref{mainA} is proved in \S\ref{mainAproof}
 and Theorem \ref{mainB} is proved in \S\ref{mainBproof}.

Before going on, we point out that in some cases these cohomology spaces are
uninteresting.
In the case $E=\C$, the cohomology groups are related, via generalizations of the Eichler--Shimura isomorphism, to certain spaces of automorphic forms.
More precisely, Franke \cite{franke} has shown that
$$
	H^{\bullet}_{\cl}(K_{f},W)
	=
	H^{\bullet}_{\relLie}(\gg,K_{\infty},W\otimes \cA(K^{f})),
$$
where $\cA(K^{f})$ is the space of automorphic forms
$\phi:\G(k)\backslash\G(\A)/K_{\infty}K_{f}\to\C$.
The right hand side is relative Lie algebra cohomology
(see for example \cite{borelwallach}).
Since the constant functions form a subspace of $\cA(K^{f})$,
 we have a $(\gg,K_{\infty})$-submodule
 $W\subset W\otimes \cA(K^{f})$.
This gives us a map:
\begin{equation}
	\label{constantcohomology}
	H^{n}_{\relLie}(\gg,K_{\infty},W)
	\to
	H^{n}_{\cl}(\G,W).
\end{equation}
We shall say that the cohomology of $\G$ is
\emph{given by constants in dimension $n$}
 if the map (\ref{constantcohomology}) is surjective.
For example the cohomology of $\SL_{2}/\Q$ is given by constants in dimensions
$0$ and $2$, although (\ref{constantcohomology}) is only bijective in dimension $0$.
On the other hand, if $\G(k)\backslash\G(\A)$ is compact then
(\ref{constantcohomology}) is injective.

It is known that the cohomology of $\G$ is given by constants
 in dimensions $n<m$ and in dimensions $n>d-m$,
 where $d$ is the common dimension of the spaces $Y(K_{f})$
 and $m$ is the real rank of $\G$.
One shows this by proving that the relative Lie algebra cohomology of
any other irreducible $(\gg,K_{\infty})$-subquotient of $W\otimes\cA(K^{f})$
vanishes in such dimensions
(see for example Corollary II.8.4 of \cite{borelwallach}).

If the cohomology is given by constants in dimension $n$, then $H^{n}_{\cl}(\G,W)$ is a finite dimensional vector space, equipped with the trivial action of $\G(\A_{f})$.
From the point of view of this paper, cohomology groups given by constants are uninteresting.
Thus Theorem \ref{crit1} is interesting only for groups of real rank $1$,
 whereas Theorem \ref{crit2} is interesting, roughly speaking, for groups of real rank $2$.

In fact we can often do a little better than Theorem \ref{variation}.
We shall prove the following in \S\ref{variation2proof}:

\begin{theorem}
	\label{variation2}
	Let $\G/k$ be connected, semi-simple and algebraically simply connected
	and assume that the weak Emerton criterion holds in dimension $n$.
	Assume also that at least one of the following two conditions holds:
	\begin{itemize}
		\item[(a)]
		$H^{p}_{\cl}(\G,\C)$ is given by constants in dimensions $p<n$
		and $H^{n+1}_{\relLie}(\gg,K_{\infty},\C)=0$; or
		\item[(b)]
		$\G(k)$ is cocompact in $\G(\A)$.
	\end{itemize}
	Then all conclusions of Theorem \ref{emertonmain}
	hold for the eigenvariety $E(n,K^{\gp})$.
\end{theorem}

The theorem is valid, for example, in the following cases
where Emerton's criterion fails:
\begin{itemize}
	\item
	$\SL_{3}/\Q$ in dimension $2$;
	\item
	$\Sp_{4}/\Q$ in dimension $2$;
	\item
	$\Spin$ groups of quadratic forms over
	 $\Q$ of signature $(2,l)$ with $l\ge 3$ in dimension $2$;
	\item
	Special unitary groups $\SU(2,l)$ with $l\ge 3$ in dimension $2$;
	\item
	$\SL_{2}/k$, where $k$ is a real quadratic field, in dimension $2$.
\end{itemize}

Our results generalize easily to simply connected, semi-simple groups as follows.
Suppose $\G/k$ is a direct sum of simply connected simple groups $\G_{i}/k$.
Assume also that the tame level $K^{\gp}$ decomposes as a direct sum of
tame levels $K^{\gp}_{i}$ in $\G_{i}(\A_{f}^{\gp})$.
By the K\"unneth formula,
 we have a decomposition of the sets of classical points:
$$
	E(n,K^{\gp})_{\cl}
	=
	\bigcup_{n_{1}+\cdots+n_{s}=n}\;
	\prod_{i=1}^{s}
	E(n_{i},K^{\gp}_{i})_{\cl}.
$$

\subsection{Some History}

Coleman and Mazur constructed the first ``eigencurve'' in \cite{colemanmazur}.
In our current notation, they constructed the $H^{1}$-eigencurve for $\GL_{2}/\Q$.
In fact they showed that the points of their eigencurve
parametrize overconvergent eigenforms.
Their arguments
were based on earlier work of Hida \cite{hida}
and Coleman \cite{coleman} on families of modular forms.
Similar results were subsequently obtained by Buzzard \cite{buzzard2004}
for the groups $\GL_{1}/k$, and for the multiplicative
group of a definite quaternion algebra over $\Q$,
 and later more generally
for totally definite quaternion algebras over totally real fields
in \cite{buzzard3}.
Kassaei \cite{kassaei-2004} treated the case that $\G$ is a form of $\GL_{2}/k$,
 where $k$ is totally real and $\G$ is split at exactly
 one archimedean place.
 Kissin and Lei in \cite{MR2154369} treated the case $\G=\GL_{2}/k$ for a totally real
 field $k$, in dimension $n=[k:\Q]$.

Ash and Stevens \cite{ashstevens, ashstevens2} obtained similar results for
$\GL_{n}/\Q$ by quite different methods.
More recently, Chenevier \cite{chenevier} constructed eigenvarieties
 for any twisted form of $\GL_{n}/\Q$
 which is compact at infinity.
Emerton's construction is apparently much more general.
However, it seems to be quite rare for his criterion to hold.
One might expect the weak criterion to hold more generally;
 in particular one might optimistically ask the following:

\begin{question}
	For $\G/k$ connected, simple, algebraically simply connected
	 and of real rank $m$, does the weak Emerton criterion always hold
	 in dimension $m$?
\end{question}

\begin{acknowledgement}
	The author benefited greatly from taking part in a
	study group on Emerton's work, organized by Kevin Buzzard.
	The author would like to thank all the participants in the London number theory
	seminar for many useful discussions.
	The author is also indebted to Prof. F. E. A. Johnson and Dr. Frank Neumann
	for their help with certain calculations.
\end{acknowledgement}

\section{Proof of Theorem \ref{crit1}}
\label{crit1proof}

Let $\G/k$ be simple, algebraically simply connected,
and assume that $\G(k_{\infty})$ is not compact.
We shall prove in \S\ref{mainAproof} that
$\tilde H^{0}(\G,E)=E$, with the trivial action of $\G(\A_{f})$.
As a consequence of this, the terms $E_{2}^{p,0}$ in
Emerton's spectral sequence (\ref{specseqlimit}) are Lie-algebra
cohomology groups of finite dimensional representations:
$$
	E_{2}^{p,0}
	=
	H^{p}_{\Lie}(\gg,W).
$$
Such cohomology groups are completely understood.
We recall some relevant results:

\begin{theorem}[Theorem 7.8.9 of \cite{weibel}]
	\label{whitehead0}
	Let $\gg$ be a semi-simple Lie algebra over a field
	of characteristic zero, and let $W$ be a finite-dimensional
	representation of $\gg$, which does not contain the
	trivial representation.
	Then we have for all $n\ge 0$,
	$$
		H^{n}_{\Lie}(\gg,W)
		=
		0.
	$$
\end{theorem}

\begin{theorem}[Whitehead's first lemma (Corollary 7.8.10 of \cite{weibel})]
	\label{whitehead1}
	Let $\gg$ be a semi-simple Lie algebra over a field
	of characteristic zero, and let $W$ be a finite-dimensional
	representation of $\gg$.
	Then we have
	$$
		H^{1}_{\Lie}(\gg,W)
		=
		0.
	$$
\end{theorem}

\begin{theorem}[Whitehead's second lemma (Corollary 7.8.12 of \cite{weibel})]
	\label{whitehead2}
	Let $\gg$ be a semi-simple Lie algebra over a field
	of characteristic zero, and let $W$ be a finite-dimensional
	representation of $\gg$.
	Then we have
	$$
		H^{2}_{\Lie}(\gg,W)
		=
		0.
	$$
\end{theorem}

We shall use these results to verify Emerton's criterion in dimension 1,
thus proving Theorem \ref{crit1}.
We must verify that the edge map \ref{edge} is an isomorphism
for $n=1$ and for every irreducible algebraic representation $W$ of $\G$.
The small terms of the spectral sequence are:
$$
	E_{2}^{\bullet,\bullet}\quad : \quad
	\begin{array}{ccc}
		\Hom_{\gg}(\check W, \tilde H^{1}(\G,E)) \medskip \\
		H^{0}_{\Lie}(\gg, W) &
		H^{1}_{\Lie}(\gg, W) &
		H^{2}_{\Lie}(\gg, W)
	\end{array}
$$
We therefore have an exact sequence:
$$
	0
	\to
	H^{1}_{\Lie}(\gg,W)
	\to
	H^{1}_{\cl}(\G,W)
	\to
	\Hom_{\gg}(\check W, \tilde H^{1}(\G,E))
	\to
	H^{2}_{\Lie}(\gg,W).
$$
By Theorems \ref{whitehead1} and \ref{whitehead2}
we know that the first and last terms are zero.
Therefore the edge map is an isomorphism.

~\sq

\section{Proof of Theorem \ref{crit2}}
\label{crit2proof}

Let $\G/k$ be connected, simple and simply connected,
and assume that $\G(k_{\infty})$ is not compact.
In \S\ref{mainBproof} we shall prove the isomorphism
$$
	\tilde{H}^{1}(\G,\Q_{p})
	=
	\Hom_{\cts}(\Cong(\G),\Q_{p})_{\G(\A_{f}^{\gp})-\smooth}.
$$
As a consequence, we have:

\begin{corollary}
	If the congruence kernel of $\G$ is finite then $\tilde H^{1}(\G,\Q_{p})=0$.
\end{corollary}

\noindent
In this context, it is worth noting that the following may be proved by a similar method.

\begin{theorem}
	If the congruence kernel of $\G$ is finite then $\tilde H^{d-1}(\G,\Q_{p})=0$,
	where $d$ is the dimension of the symmetric space $\G(k_{\infty})/K_{\infty}$.
\end{theorem}

We shall use the corollary to verify the weak Emerton criterion in
dimension $2$.
Suppose first that $W$ is a non-trivial irreducible
algebraic representation of $\G$. We must show that the edge map
(\ref{edge}) is an isomorphism.
By Theorem \ref{whitehead0} we know that the bottom row of the spectral sequence
is zero, and by the corollary we know that the first row is zero.
The small terms of the spectral sequence are as follows:
$$
	E_{2}^{\bullet,\bullet}\quad :\quad
	\begin{array}{cccc}
		\Hom_{\gg}(\check W, \tilde H^{2}(\G,E)_{\la})
		\medskip\\
		0&0&0
		\medskip\\
		0&0&0&0
	\end{array}
$$
Hence in this case the edge map is an isomorphism.

In the case that $W$ is the trivial representation,
we must only verify that the edge map is injective and that its cokernel is
a finite dimensional trivial representation of $\G(\A_{f})$.
We still know in this case that the first row of the spectral sequence is zero.
For the bottom row, Theorems \ref{whitehead1} and \ref{whitehead2}
tell us that the spectral sequence is as follows:
$$
	E_{2}^{\bullet,\bullet}\quad :\quad
	\begin{array}{cccc}
		\Hom_{\gg}(E, \tilde H^{2}(\G,E)_{\la})
		\medskip\\
		0&0&0
		\medskip\\
		E&0&0&H^{3}_{\Lie}(\gg,E)
	\end{array}
$$
It follows that we have an exact sequence
\begin{equation}
	\label{enigma}
	0
	\to
	H^{2}_{\cl}(\G,E)
	\to
	\Hom_{\gg}(E, \tilde H^{2}(\G,E)_{\la})
	\to
	H^{3}_{\Lie}(\gg,E).
\end{equation}
The action of $\G(\A_{f})$ on $H^{3}_{\Lie}(\gg,E)$ is trivial, since this action is
defined by the (trivial) action on $\tilde H^{0}(\G,E)=E$.

~\sq
\medskip

\begin{remark}
	It is interesting to calculate the cokernel of the edge map
	in (\ref{enigma}).
	In fact it is known that for any simple Lie algebra $\gg$ over a
	field $E$ of characteristic zero, $H^{3}_{\Lie}(\gg,E)=E$.
	We therefore have by the K\"unneth formula:
	$$
		H^{3}_{\Lie}(\gg, E)
		=
		E^{d},
	$$
	where $d$ is the number of simple factors of
	$\G\times_{k} \bar k$.
	In particular, this is never zero.
	The exact sequence (\ref{enigma})
	can be continued for another term as follows:
	$$
		0
		\to
		H^{2}_{\cl}(\G,E)
		\to
		\tilde{H}^{2}(\G,E)_{\la}^{\gg}
		\to
		H^{3}_{\Lie}(\gg,E)
		\to
		H^{3}_{\cl}(\G,E)^{\G(\A_{f})}.
	$$
	In order to calculate the last term, we first choose an embedding of $E$
	in $\C$, and tensor with $\C$.
	There is a map
	$$
		H^{3}_{\relLie}(\gg,K_{\infty},\C)
		\to
		H^{3}_{\cl}(\G,\C)^{\G(\A_{f})}.
	$$
	If the $k$-rank of $\G$ is zero, then this map is an isomorphism.
	In other cases, it is often surjective,
	although the author does not know how to prove this statement in general.
	The groups $H^{\bullet}_{\relLie}(\gg, K_{\infty},\C)$ are the cohomology groups of
	compact symmetric spaces (see \S I.1.6 of \cite{borelwallach})
	and are completely understood.
	In particular, it is often the case that $H^{3}_{\relLie}(\gg,K_{\infty},\C)=0$.
	This implies that the edge map in (\ref{enigma})
	often has a non-trivial cokernel.
\end{remark}

\section{Proof of Theorem \ref{variation}}
\label{variationproof}

Theorem \ref{variation} is a variation on Theorem \ref{emertonmain}.
In order to prove it, we recall some of the intermediate steps in Emerton's
proof of Theorem \ref{emertonmain}.

In \cite{emertonJacq}, Emerton introduced a new kind of Jacquet functor,
$\Jacq_{\B}$, from the category of essentially admissible
(in the sense of Definition 6.4.9 of \cite{emertonAnalytic})
locally analytic representations of $\G(k_{\gp})$
to the category of essentially admissible
locally analytic representations of $\T(k_{\gp})$.
This functor is left exact, and its restriction to the full subcategory
 of smooth representations is exact.
Indeed, its restriction to smooth representations is the usual
Jacquet functor of coinvariants.

Applying the Jacquet functor to the space $\tilde H^{n}(K^{\gp},E)_{\la}$,
one obtains an essentially admissible locally analytic representation of
$\T(k_{\gp})$.
On the other hand, the category of essentially admissible
 locally analytic representations of $\T(k_{\gp})$
 is anti-equivalent to the category of
 coherent rigid analytic sheaves on $\hat T$
 (Proposition 2.3.2 of \cite{emerton}).
We therefore have a coherent sheaf $\scE$ on $\hat T$.
Since the action of $\cH(K^{\gp})$ on $\tilde H^{n}(K^{\gp},E)_{\la}$ commutes with
 that of $\G(k_{\gp})$, it follows that $\cH(K^{\gp})$
 acts on $\scE$.
Let $\scA$ be the image of $\cH(K^{\gp})^{\sph}$
 in the sheaf of endomorphisms of $\scE$.
Thus $\scA$ is a coherent sheaf of commutative rings on $\hat T$.
Writing $\Spec\scA$ for the relative spec of $\scA$ over $\hat T$,
 we have a Zariski-closed embedding
 $\Spec\scA\to \hat T \times \Spec\cH(K^{\gp})^{\sph}$.
Since $\scA$ acts as endomorphisms of $\scE$,
 we may localize $\scE$ to a coherent sheaf $\scM$ on $\Spec \scA$.

Theorem \ref{emertonmain} may be deduced from the following two results.

\begin{theorem}[2.3.3 of \cite{emerton}]
	\label{emA}
	\begin{itemize}
		\item[(i)]
		The natural projection $\Spec \scA \to \hat T$ is a finite morphism. 
		\item[(ii)]
		The map $\Spec \scA \to \check \gt$ has discrete fibres.
		\item[(iii)]
		The fibre of $\scM$ over a point $(\chi,\lambda)$
		of $\hat T \times \Spec \cH ( K^{\gp})^{\sph}$ is dual 
		to the $(\T(k_{\gp}) = \chi, \cH ( K^{\gp})^{\sph}= \lambda)$-eigenspace
		of $\Jacq_{\B} ( \tilde H^{n}( K^{\gp} ,E)_{\la})$.
		In particular, the point $(\chi, \lambda)$
		lies in $\Spec \scA$ if and only if this eigenspace is non-zero.
	\end{itemize}
\end{theorem}

For any representation $V$ of $\G(k_{\gp})$ over $E$,
we shall write $V_{W-\lalg}$ for the subspace of $W$-locally
algebraic vectors in $V$.
Note that under Emerton's criterion, we have
\begin{equation}
	\label{Walgebraicvectors}
	H^{n}_{\cl}(K^{\gp},W)\otimes \check W
	=
	\tilde H^{n}(K^{\gp},E)_{\check W-\lalg}.
\end{equation}
Hence $H^{n}_{\cl}(K^{\gp},W)\otimes \check W$ is a closed subspace
of $\tilde H^{n}(K^{\gp},E)_{\la}$.
By left-exactness of $\Jacq_{\B}$ we have an injective map
$$
	\Jacq_{\B}(H^{n}_{\cl}(K^{\gp},W)\otimes \check W)
	\to
	\Jacq_{\B}(\tilde H^{n}(K^{\gp},E)_{\la})
$$
There are actions of $\T(k_{\gp})$ and $\cH(K^{\gp})$ on these
spaces, so we may restrict this map to eigenspaces:
$$
	\Jacq_{\B}(H^{n}_{\cl}(K^{\gp},W)\otimes \check W)^{(\chi,\lambda)}
	\to
	\Jacq_{\B}(\tilde H^{n}(K^{\gp},E)_{\la})^{(\chi,\lambda)},
	\quad
	(\chi,\lambda)\in \hat T \times \Spec\cH(K^{\gp})^{\sph}.
$$
The next result tells us that this restriction is often an isomorphism.

\begin{theorem}[Theorem 4.4.5 of \cite{emertonJacq}]
	\label{emB}
	Let $V$ be an admissible continuous representation
	of $\G(k_{\gp})$ on a Banach space.
	If $\chi := \theta\psi_{W}\in \hat T(\bar E)$ is of non-critical slope,
	then the closed embedding
	$$
		\Jacq_{\B}(V_{W-\lalg} )
		\to
		\Jacq_{\B}(V_{\la})
	$$
	induces an isomorphism on $\chi$-eigenspaces.
\end{theorem}

We recall Theorem \ref{variation}.

\begin{thm}
	If the weak Emerton criterion holds for $\G$ in dimension $n$,
	then
	\begin{enumerate}
		\item
		Projection onto the first factor induces a finite map
		$E(n,K^{\gp})\to\hat{T}$.
		\item
		The map $E(n,K^{\gp})\to\check{\gt}$ has discrete fibres.
		\item
		If $(\chi,\lambda)$ is a point of the Eigencurve such that $\chi$ is
		locally algebraic and of non-critical slope,
		then either $(\chi,\lambda)$ is a classical point
		or $(\chi,\lambda)$ is the trivial representation
		of $\T(k_{\gp})\times\cH(K^{\gp})^{\sph}$.
	\end{enumerate}
\end{thm}

\proof
To prove the first two parts of the theorem, it is sufficient to
show that $E(K^{\gp},n)$ is a closed subspace of $\Spec\scA$.
Since $E(n, K^{\gp})$ is defined to be the closure of the set of
 classical points, it suffices to show that each classical point
 is in $\Spec\scA$.

Suppose $\pi$ is a representation appearing in $H^{n}_{\cl}(K^{\gp},W)$
and let $(\theta\psi_{W}, \lambda)$ be the corresponding classical point.
This means that the $(\theta,\lambda)$-eigenspace in the $\Jacq_{\B}(\pi)$
is non-zero.
By exactness of the Jacquet functor on smooth representations,
 it follows that the $(\theta,\lambda)$ eigenspace in the Jacquet module
of $H^{n}_{\cl}(K^{\gp},W)$ is non-zero.
Hence by Proposition 4.3.6 of \cite{emertonJacq},
the $(\theta\psi_{W},\lambda)$-eigenspace in the Jacquet module of
$H^{n}_{\cl}(K^{\gp},W)\otimes \check W$
is non-zero.
By left-exactness of the Jacquet functor, it follows that the
$(\theta\psi_{W},\lambda)$ eigenspace in the Jacquet module of
$\tilde H^{n}(\G,E)_{\la}$ is non-zero.
Hence by Theorem \ref{emA} (iii) it follows that
the classical point is in $\Spec\scA$.

If $(\theta\psi,\lambda)$ is of non-critical slope
then Theorem \ref{emB} shows that the converse also holds.
\sq
\medskip

\section{Proof of Theorem \ref{variation2}}
\label{variation2proof}

We first recall the statement:

\begin{thm}
	Let $\G/k$ be connected, semi-simple and algebraically simply connected
	and assume that the weak Emerton criterion holds in dimension $n$.
	Assume also, that at least one of the following two conditions holds:
	\begin{itemize}
		\item[(a)]
		$H^{p}_{\cl}(\G,\C)$ is given by constants in dimensions $p<n$
		and $H^{n+1}_{\relLie}(\gg,K_{\infty},\C)=0$; or
		\item[(b)]
		$\G(k)$ is cocompact in $\G(\A)$.
	\end{itemize}
	Then all the conclusions of Theorem \ref{emertonmain}
	hold for the eigenvariety $E(n,K^{\gp})$.
\end{thm}

\proof
To prove the theorem, we shall find a
 continuous admissible Banach space representation
 $V$,
 such that for every irreducible algebraic representation $W$,
 there is an isomorphism of smooth $\G(\A_{f})$-modules
\begin{equation}
	\label{toprove}
	H^{n}_{\cl}(\G,W)
	\cong
	\Hom_{\gg}(\check W,V_{\la}).
\end{equation}
Recall that by the weak Emerton criterion,
 we have an exact sequence of smooth $\G(\A_{f})$-modules
\begin{equation}
	\label{weaksequ1}
	0
	\to
	H^{n}_{\cl}(\G,E)
	\to
	\tilde H^{n}(\G,E)_{\la}^{\gg}
	\to
	E^{r}
	\to
	0,
	\qquad
	r\ge 0.
\end{equation}
It follows, either from Lemma \ref{keylemma} or from Lemma
 \ref{keylemma2} below,
 that all such sequences split.
We therefore have a subspace $E^{r}\subset \tilde H^{n}(\G,E)$,
on which $\G(\A_{f})$ acts trivially.
We define $V$ to be the quotient, so that there is an exact sequence
of admissible continuous representations of $\G(\A_{f})$ on $E$-Banach spaces.
\begin{equation}
	\label{ex1}
	0
	\to
	E^{r}
	\to
	\tilde H^{n}(\G,E)
	\to
	V
	\to
	0.
\end{equation}
Taking $\gg$-invariants of (\ref{ex1})
 and applying Whitehead's first lemma (Theorem \ref{whitehead1}),
 we have an exact sequence:
\begin{equation}
	\label{weaksequ2}
	0
	\to
	E^{r}
	\to
	\tilde H^{n}(\G,E)_{\la}^{\gg}
	\to
	V_{\la}^{\gg}
	\to
	0.
\end{equation}
On the other hand, $E^{r}$ is a direct summand of $\tilde H^{n}(\G,E)^{\gg}_{\la}$,
so this sequence also splits.
Comparing (\ref{weaksequ1}) and (\ref{weaksequ2}),
we obtain 
$$
	H^{n}_{\cl}(\G,E)
	=
	V_{\la}^{\gg}
	=
	\Hom_{\gg}(E,V_{\la}).
$$
This verifies (\ref{toprove}) in the case that $W$ is the trivial representation.

Now taking $W$ to be a non-trivial irreducible representation, and applying
$\Hom_{\gg}(\check W,-_{\la})$ to (\ref{ex1}), we obtain a long exact sequence:
$$
	0
	\to
	\Hom_{\gg}(\check W,\tilde H^{n}(\G,E)_{\la})
	\to
	\Hom_{\gg}(\check W,V_{\la})
	\to
	\Ext_{\gg}^{1}(\check W,E^{r}).
$$
By Whitehead's first lemma, the final term above is zero.
Hence, by the weak Emerton criterion, we have:
$$
	H^{2}_{\cl}(\G,W)
	=
	\Hom_{\gg}(\check W,\tilde H^{n}(\G,E)_{\la})
	=
	\Hom_{\gg}(\check W,V_{\la}).
$$
\sq
\medskip

\begin{lemma}
	\label{keylemma}
	Assume that $H^{q}_{\cl}(\G,\C)$ is given by constants in
	dimensions $q<n$ and $H_{\relLie}^{n+1}(\gg,K_{\infty},\C)=0$.
	Then
	$$
		\Ext^{1}_{\G(\A_{f})}(E, H^{n}_{\cl}(\G,E))=0,
	$$
	where the $\Ext$-group is calculated from the category of
	smooth representations of $\G(\A_{f})$ over $E$.
\end{lemma}

\proof
Since we are dealing with smooth representations, the topology of $E$ plays no role, so it is sufficient to prove that
$$
	\Ext^{1}_{\G(\A_{f})}(\C, H^{n}_{\cl}(\G,\C))=0,
$$
To prove this, it is sufficient to show that for every sufficiently large
 finite set $S$ of finite primes of $k$, we have
$$
	\Ext^{1}_{\G(k_{S})}(\C, H^{n}_{\cl}(\G,\C))=0.
$$
For this, we shall use the spectral sequence of Borel (\S3.9 of \cite{borel2};
see also \S2 of \cite{blasiusfrankegrunewald}):
$$
	\Ext^{p}_{\G(k_{S})}(E, H^{q}_{\cl}(\G,\C))
	\implies
	H^{p+q}_{S-\cl}(\G,\C),
$$
where $H^{\bullet}_{S-\cl}(\G,-)$ denotes the direct limit over
all $S$-congruence subgroups:
$$
	H^{\bullet}_{S-\cl}(\G,-)
	=
	\lim_{\to\atop K^{S}}
	H^{\bullet}_{\Group}(\Gamma_{S}(K^{S}),-),
	\qquad
	\Gamma_{S}(K_{S})
	=
	\G(k)\cap \left(\G(k_{\infty\cup S})\times K^{S}\right).
$$
By Proposition X.4.7 of \cite{borelwallach}, we know that
$$
	\Ext^{p}_{\G(k_{S})}(\C,\C)=0,
	\qquad
	p\ge 1.
$$
Since $H^{q}_{\cl}(\G,\C)$ is a trivial representation of $\G(\A_{f})$ in
dimensions $q<n$, it follows from Borel's spectral sequence
that $\Ext^{1}_{\G(k_{S})}(\C, H^{n}_{\cl}(\G,\C))$ injects into $H^{n+1}_{S-\cl}(\G,\C)$.
On the other hand,
 it is shown in Theorem 1 of \cite{blasiusfrankegrunewald},
 that for $S$ sufficiently large, $H^{n+1}_{S-\cl}(\G,\C)$ is isomorphic
 to $H^{n+1}_{\relLie}(\gg,K_{\infty},\C)$.
 
Under the hypothesis that $H^{n+1}_{\relLie}(\gg,K_{\infty},\C)=0$,
 it follows that for $S$ sufficiently large,
 $\Ext^{1}_{\G(k_{S})}(\C, H^{n}_{\cl}(\G,\C))=0$.
\sq
\medskip

\begin{lemma}
	\label{keylemma2}
	Assume that $\G(k)$ is cocompact in $\G(\A)$.
	Then
	$$
		\Ext^{1}_{\G(\A_{f})}(E, H^{n}_{\cl}(\G,E))=0,
	$$
	where the $\Ext$-group is calculated from the category of
	smooth representations of $\G(\A_{f})$ over $E$.
\end{lemma}

(The argument in fact shows that $\Ext^{p}_{\G(\A_{f})}(E, H^{q}_{\cl}(\G,E))=0$
for all $p>0$.)
\medskip

\proof
As in the proof of the previous lemma,
 we shall show that for $S$ sufficiently large,
$$
	\Ext^{1}_{\G(k_{S})}(\C,H^{n}_{\cl}(\G,\C))=0.
$$
Recall that we have a decomposition:
$$
	L^{2}(\G(k)\backslash\G(\A))
	=
	\widehat
	{\bigoplus_{\pi}}\
	m(\pi) \cdot \pi,
$$
with finite multiplicities $m(\pi)$ and automorphic representations $\pi$.
Here the $\hat\oplus$ denotes a Hilbert space direct sum.
We shall write $\pi=\pi_{\infty}\otimes\pi_{f}$,
 where $\pi_{\infty}$ is an irreducible unitary representation
 of $\G(k_{\infty})$ and $\pi_{f}$ is a smooth irreducible unitary
 representation of $\G(\A_{f})$.
This decomposition may by used to calculate the classical cohomology
(Theorem VII.6.1 of \cite{borelwallach}):
$$
	H^{\bullet}_{\cl}(\G,\C)
	=
	\bigoplus_{\pi}
	m(\pi) \cdot H^{\bullet}_{\relLie}(\gg,K_{\infty},\pi_{\infty})\otimes \pi_{f}.
$$
It is therefore sufficient to show that for each automorphic representation $\pi$,
we have (for $S$ sufficiently large) $\Ext^{1}_{\G(k_{S})}(\C,\pi_{f})=0$.
The smooth representation $\pi_{f}$ decomposes as a tensor product
of representations of $\G(k_{\gq})$ for $\gq\in S$, together with a representation
of $\G(\A_{f}^{S})$:
$$
	\pi_{f}
	=
	\left(\bigotimes_{\gq\in S}\pi_{\gq}\right)
	\otimes
	\pi_{f}^{S}.
$$
This gives a decomposition of the cohomology:
\begin{equation}
	\label{automorphictensorprod}
	\Ext^{\bullet}_{\G(k_{S})}(\C,\pi_{f})
	=
	\left(\bigotimes_{\gq\in S} \Ext^{\bullet}_{\G(k_{\gq})}(\C,\pi_{\gq})\right)
	\otimes
	\pi_{f}^{S}.
\end{equation}
There are two cases to consider.

\emph{Case 1.}
Suppose $\pi$ is the trivial representation,
 consisting of the constant functions on $\G(k)\backslash \G(\A)$.
Then by Proposition X.4.7 of \cite{borelwallach},
 we have
$$
	\Ext^{n}_{\G(k_{\gq})}(\C,\C)
	=
	0,
	\qquad
	n\ge 1.
$$
This implies by (\ref{automorphictensorprod}) that $\Ext^{1}_{\G(k_{S})}(\C,\C)=0$.

\emph{Case 2.}
Suppose $\pi$ is non-trivial, and hence contains no non-zero constant functions.
If $\gq$ is a prime for which no factor of $\G(k_{\gq})$ is compact,
 then it follows from the strong approximation theorem that
 the local representation $\pi_{\gq}$ is non-trivial.
This implies that
$$
	\Ext^{0}_{\G(k_{\gq})}(\C,\pi_{\gq})
	=
	\Hom_{\G(k_{\gq})}(\C,\pi_{\gq})
	=
	0.
$$
If $S$ contains at least two such primes, then we have
 by (\ref{automorphictensorprod})
$$
	\Ext^{1}_{\G(k_{S})}(\C,\pi_{f})
	=
	0.
$$
\sq
\medskip

\begin{remark}
	At first sight, it might appear that $\Ext^{1}_{\G(\A_{f})}(\C, H^{n}_{\cl}(\G,\C))$
	 should always be zero; however this is not the case.
	For example, if $\G=\SL_{2}/\Q$ then
	$$
		\Ext^{1}_{\SL_{2}(\A_{f})}(\C,H^{1}_{\cl}(\SL_{2}/\Q,\C))
		=
		\C.
	$$
	This may be verified using the spectral sequence of Borel cited above,
	together with the fact that $H^{2}_{\relLie}(\sll_{2},\SO(2),\C)=\C$.
\end{remark}

\section{Proof of Theorem \ref{mainA}}
\label{mainAproof}

We assume in this section that $\G/k$ is connected, simple
and algebraically simply connected, and that $\G(k_{\infty})$ is not compact.

\begin{proposition}
	As topological spaces, we have
	$Y(K_{f})=\Gamma(K_{f}) \backslash \G(k_{\infty})/K_{\infty}$.
\end{proposition}

\proof
By the strong approximation theorem
(Theorem \ref{SAT}),
$\G(k)\G(k_{\infty})$ is a dense subgroup of $\G(\A)$.
Since $\G(k_{\infty})K_{f}$ is open in $\G(\A)$, this implies
that $\G(k)\G(k_{\infty})K_{f}$ is a dense, open subgroup
of $\G(\A)$.
Since open subgroups are closed it follows that
$$
	\G(k)\G(k_{\infty})K_{f}
	=
	\G(\A).
$$
Quotienting out on the left by $\G(k)$,
we have (as coset spaces):
$$
	(\G(k)\cap \G(k_{\infty})K_{f}) \backslash ( \G(k_{\infty})K_{f} )
	=
	\G(k) \backslash \G(\A).
$$
Substituting the definition of $\Gamma(K_{f})$, we have:
$$
	\Gamma(K_{f}) \backslash \G(k_{\infty})K_{f}
	=
	\G(k)\backslash\G(\A).
$$
Quotienting out on the right by $K_{\infty}K_{f}$,
 we get:
$$
	\Gamma(K_{f}) \backslash \G(k_{\infty}) /K_{\infty}
	=
	Y(K_{f}).
$$
\sq
\medskip

In particular, this implies:

\begin{corollary}
	$Y(K_{f})$ is connected.
\end{corollary}

\proof
$\G(k_{\infty})$ is connected.
\sq
\medskip

If $K_{f}$ is sufficiently small then the group $\Gamma(K_{f})$ is torsion-free.
We shall assume that this is the case.
Hence $Y(K_{f})$ is a manifold.
Its universal cover is $\G(\R)/K$,
and its fundamental group is $\Gamma(K_{f})$.

\begin{corollary}
	\label{classgroup}
	If $\Gamma(K_{f})$ is torsion-free then
	$H^{\bullet}(Y(K_{f}), - )=H^{\bullet}_{\Group}(\Gamma(K_{f}), -)$.
\end{corollary}

\proof
This follows because $\Gamma(K_{f})$ is the fundamental group of
$Y(K_{f})$, and the universal cover $\G(k_{\infty})/K_{\infty}$
is contractible.
See for example \cite{serre2}.
\sq
\medskip

\begin{corollary}
	Let $\G/k$ be connected, simple, simply connected
	 and assume $\G(k_{\infty})$ is not compact.
	Then as $\G(\A_{f})$-modules,
	$$
		\tilde H^{0}(\G,E)_{\la}
		=
		\tilde H^{0}(\G,E)
		=
		E.
	$$
\end{corollary}

\proof
Since every $Y(K_{f})$ is connected, we have a canonical isomorphism:
$$
	H^{0}(Y(K_{\gp}K^{\gp}),\Z/p^{s})
	=
	\Z/p^{s}.
$$
Furthermore, the pull-back maps
$$
	H^{0}(Y(K_{\gp}K^{\gp}),\Z/p^{s})
	\to
	H^{0}(Y(K_{\gp}'K^{\gp}),\Z/p^{s})
	\qquad
	(K_{\gp}'
	\subset
	K_{\gp})
$$
are all the identity on $\Z/p^{s}$.
It follows that
$$
	\limdKp
	H^{0}(Y(K^{\gp}K_{\gp}),\Z/p^{s})
	=
	\Z/p^{s}.
$$
Since the pull-back maps are all the identity, it follows
that the action of $\G(k_{\gp})$ on this group is trivial.
Taking the projective limit over $s$ and tensoring with $E$ we
find that
$$
	\tilde H^{0}(K^{\gp},E)
	=
	E.
$$
The action of $\G(k_{\gp})$ is clearly still trivial,
and hence every vector is locally analytic.
The groups $\tilde H^{0}(K^{\gp},E)$ for varying tame level $K^{\gp}$
 form a direct system with respect to the pullback maps.
These pullback maps are all the identity on $E$.
Taking the direct limit over the tame levels, we obtain:
$$
	\tilde H^{0}(\G,E)
	=
	E.
$$
Since the pullback maps are all the identity,
it follows that the action of $\G(\A_{f}^{\gp})$ on $\tilde H^{0}(\G,E)$ is trivial.
\sq
\medskip

\section{Some Cohomology Theories}

In this section we introduce some notation and recall some results,
 which will be needed in the proof of Theorem \ref{mainB}.

\subsection{Discrete cohomology}

Let $G$ be a profinite group acting on an abelian group $A$.
We say that the action is \emph{smooth} if every element of $A$ has
open stabilizer in $G$.
For a smooth $G$-module $A$, we define
$H^{\bullet}_{\disc}(G,A)$ to be the cohomology
 of the complex of smooth cochains on $G$ with values in $A$.
Due to compactness,
 cochains take only finitely many values, so we have
$$
	H^{\bullet}_{\disc}(G,A)
	=
	\lim_{\rightarrow\atop U}
	H^{\bullet}_{\Group}(G/U,A^{U}).
$$
Here the limit is taken over the open normal subgroups $U$ of $G$,
and the cohomology groups on the right hand side are those of finite groups.

\begin{theorem}[Hochschild--Serre spectral sequence
	(\S2.6b of \cite{serre})]
	\label{hochschildserre}
	Let $G$ be a profinite group and $A$ a discrete $G$-module
	on which $G$ acts smoothly.
	Let $H$ be a closed, normal subgroup.
	Then there is a spectral sequence:
	$$
		H^{p}_{\disc}(G/H,H^{q}_{\disc}(H,A))
		\implies
		H^{p+q}_{\disc}(G,A).
	$$
\end{theorem}

For calculations with ad\`ele groups, we need the following result on
countable products of groups.

\begin{proposition}[see \S2.2 of \cite{serre}]
	Let
	$$
		G
		=
		\prod_{i\in\N}
		G_{i}
	$$
	be a countable product of profinite groups
	and let $A$ be a discrete $G$-module.
	For any finite subset $S\subset \N$ we let
	$$
		G_{S}
		=
		\prod_{i\in S}
		G_{i}.
	$$
	Then
	$$
		H^{n}_{\disc}(G,A)
		=
		\lim_{\to \atop S}
		H^{n}_{\disc}(G_{S},A).
	$$
	Here the limit is taken over all finite subsets with respect to the inflation maps.
\end{proposition}

\begin{corollary}
	\label{localglobal}
	Let $G$ and $A$ be as in the previous
	proposition,
	and assume that the action of $G$ on $A$ is trivial.
	Assume also that for a fixed $n$, we have:
	$$
		H^{r}_{\disc}(G_{i},A)
		=
		0,
		\quad
		r=1,\ldots,n-1,\;
		i\in\N.
	$$
	Then
	$$
		H^{n}_{\disc}(G,A)
		=
		\bigoplus_{i\in\N}
		H^{n}_{\disc}(G_{i},A).
	$$
\end{corollary}

\proof
Let $S\subset\N$ be a finite set and let $i\notin S$.
We have a direct sum decomposition
$$
	G_{S\cup \{i\}}
	=
	G_{S}\oplus G_{i}.
$$
Regarding this as a (trivial) group extension,
we have a spectral sequence:
$$
	H^{p}_{\disc}(G_{S},H^{q}(G_{i},A))
	\implies
	H^{p+q}_{\disc}(G_{S\cup\{i\}},A).
$$
since the sum is direct, it follows that all the maps in the spectral sequence are zero,
and we have
$$
	H^{n}_{\disc}(G_{S\cup\{i\}},A)
	=
	\bigoplus_{r=0}^{n}
	H^{n-r}_{\disc}(G_{S},H^{r}_{\disc}(G_{i},A)).
$$
By our hypothesis, most of these terms vanish, and we are left with:
$$
	H^{n}_{\disc}(G_{S\cup\{i\}},A)
	=
	H^{n}_{\disc}(G_{S},A)\oplus H^{n}_{\disc}(G_{i},A).
$$
By induction on the size of $S$, we deduce that
$$
	H^{n}_{\disc}(G_{S},A)
	=
	\bigoplus_{i\in S}
	H^{n}_{\disc}(G_{i},A).
$$
The result follows from the previous proposition.
\sq
\medskip

\subsection{Continuous cohomology}

Again suppose that $G$ is a profinite group, acting on an abelian topological
group $A$.
We call $A$ a continuous $G$-module if the map $G\times A \to A$ is
continuous.
For a continuous $G$-module $A$,
 we define the continuous cohomology groups $H^{\bullet}_{\cts}(G,A)$
 to be the cohomology of the complex of continuous cochains.
If the topology on $A$ is actually discrete then continuous cochains
are in fact smooth, so we have
$$
	H^{\bullet}_{\cts}(G,A)
	=
	H^{\bullet}_{\disc}(G,A).
$$

\subsection{Derived functors of inverse limit}

Let $\Ab$ be the category of abelian groups.
By a projective system in $\Ab$, we shall mean a collection of
objects $A_{s}$ ($s\in\N$) and morphisms $\phi:A_{s+1}\to A_{s}$.
We shall write $\Ab^{\N}$ for the category of projective systems in $\Ab$.
There is a functor
$$
	\limprojs:\Ab^{\N}\to \Ab.
$$
This functor is left-exact. It has right derived functors
$$
	\left(\limprojs\right)^{\bullet}:\Ab^{\N}\to \Ab.
$$
It turns out that $\left(\limprojs\right)^{n}$ is zero for $n\ge 2$.
The first derived functor has the following simple description due to Eilenberg.
We define a homomorphism
$$
	\Delta:
	\prod_{s}A_{s}
	\to
	\prod_{s}A_{s},
	\qquad
	\big(\Delta(a_{\bullet})\big)_{s}
	=
	a_{s}-\phi(a_{s+1}).
$$
With this notation we have
$$
	\limprojs\ A_{s}
	=
	\ker\Delta.
$$
Eilenberg showed that
$$
	\limprojsd A_{s}
	=
	\coker\Delta.
$$
A projective system $A_{s}$ is said to satisfy the \emph{Mittag--Leffler condition}
if for every $s\in \N$ there is a $t\ge s$ such that for every $u\ge t$
the image of $A_{u}$ in $A_{s}$ is equal to the image of $A_{t}$ in $A_{s}$.

\begin{proposition}[Proposition 3.5.7 of \cite{weibel}]
	If $A_{s}$ satisfies the Mittag--Leffler condition then
	$\limprojsd A_{s}=0$.
\end{proposition}

This immediately implies:

\begin{corollary}[Exercise 3.5.2 of \cite{weibel}]
	\label{finitegroups}
	If $A_{s}$ is a projective system of finite abelian groups
	 then	$\limprojsd A_{s}=0$.
\end{corollary}

We shall use the derived functor $\limprojsd$ to pass
between discrete and continuous cohomology:

\begin{theorem}[Eilenberg--Moore Sequence (Theorem 2.3.4 of \cite{neukirch})]
	\label{eilenbergmoorethm}
	Let $G$ be a profinite group and $A$ a projective limit
	of finite discrete continuous $G$-modules
	$$
		A
		=
		\limprojs A_{s}.
	$$
	Then there is an exact sequence:
	$$
		0
		\to
		\limprojsd 
		H^{n-1}_{\disc}(G,A_{s})
		\to
		H^{n}_{\cts}(G,A)
		\to
		\limprojs
		H^{n}_{\disc}(G,A_{s})
		\to
		0.
	$$
\end{theorem}

\subsection{Stable Cohomology}

For a continuous representation $V$ of $\G(k_{\gp})$ over $E$,
 we shall write $V_{\st}$ for the set of smooth vectors.
The functor $V\mapsto V_{\st}$ is left exact from the category of
 continuous admissible representations of $\G(k_{\gp})$
 (in the sense of \cite{schneider-teitelbaum})
 to the category of smooth representations.
We shall write $H^{\bullet}_{\st}(\G(k_{\gp}),-)$ for the right-derived
functors.
This is called ``stable cohomology'' by Emerton (Definition 1.1.5 of \cite{emerton}).
It turns out that stable cohomology may be expressed in terms
of continuous group cohomology as follows (Proposition 1.1.6 of \cite{emerton}):
$$
	H^{\bullet}_{\st}(\G(k_{\gp}),V)
	=
	\limdKp
	H^{\bullet}_{\cts}(K_{\gp},V).
$$
There is an alternative description of these derived functors which we shall also use.
Let $V_{\la}$ denote the subspace of locally analytic vectors in $V$.
There is an action of the Lie algebra $\gg$ on
 $V_{\la}$.
Stable cohomology may be expressed in terms of Lie
algebra cohomology as follows
(Theorem 1.1.13 of \cite{emerton}):
\begin{equation}
	\label{stableLie}
	H^{\bullet}_{\st}(\G(k_{\gp}),V)
	=
	H^{\bullet}_{\Lie}(\gg,V_{\la}).
\end{equation}

\section{Proof of Theorem \ref{mainB}}
\label{mainBproof}

In this section, we shall assume that $\G/k$ is connected, simply connected and simple, and that $\G(k_{\infty})$ is not connected.
We regard the vector space $\Hom_{\cts}(\Cong(\G),\Q_{p})$
as a $p$-adic Banach space with the supremum norm:
$$
	|| \phi ||
	=
	\sup_{x\in \Cong(\G)} |\phi(x)|_{p}.
$$
We regard $\Hom_{\cts}(\Cong(\G),\Z/p^{s})$ as a discrete
abelian group.
The group $\G(\A_{f})$ acts on these spaces as follows:
$$
	(g\phi)(x)
	=
	\phi(g^{-1}x g),
	\quad
	g\in \G(\A_{f}),\; x\in \Cong(\G).
$$

\begin{lemma}
	\label{smooth}
	The action of $\G(\A_{f})$ on $\Hom_{\cts}(\Cong(\G),\Z/p^{s})$ is smooth.
\end{lemma}

\proof
One may prove this directly; however it is implicit in the Hochschild--Serre
spectral sequence.
It is sufficient to show that the action of some open subgroup is smooth.
Let $K_{f}$ be a compact open subgroup of $\G(\A_{f})$, and write
write $\tilde K_{f}$ for the preimage of $K_{f}$ in $\tilde \G(k)$.
We therefore have an extension of profinite groups:
$$
	1
	\to
	\Cong(\G)
	\to
	\tilde K_{f}
	\to
	K_{f}
	\to
	1.
$$
We shall regard $\Z/p^{s}$ as a trivial, and hence smooth, $\tilde K_{f}$-module.
It follows that each $H^{q}_{\disc}(\Cong(\G),\Z/p^{s})$ is a smooth $K_{f}$-module.
On the other hand we have
$$
	\Hom_{\cts}(\Cong(\G),\Z/p^{s})
	=
	H^{1}_{\disc}(\Cong(\G),\Z/p^{s}).
$$
\sq

\begin{lemma}
	The action of $\G(\A_{f})$ on the $p$-adic Banach space
	 $\Hom_{\cts}(\Cong(\G),\Q_{p})$ is continuous.
\end{lemma}

\proof
It is sufficient to prove this for the open submodule
$\Hom_{\cts}(\Cong(\G),\Z_{p})$.
We have, as topological $\G(\A_{f})$-modules:
$$
	\Hom_{\cts}(\Cong(\G),\Z_{p})
	=
	\limprojs
	\Hom_{\cts}(\Cong(\G),\Z/p^{s}).
$$
Continuity follows from the previous Lemma.
\sq
\medskip

We shall say that a vector $v\in \Hom_{\cts}(\Cong(\G),\Q_{p})$ is
$\G(\A_{f}^{\gp})$-smooth if its stabilizer in $\G(\A_{f}^{\gp})$ is open.
The set of such vectors will be written
$$
	\Hom_{\cts}(\Cong(\G),\Q_{p})_{\G(\A_{f}^{\gp})-\smooth}.
$$

\begin{thm}
	Assume $\G/k$ is connected, simple and simply connected, and
	that $\G(k_{\infty})$ is not compact.
	Then we have an isomorphism of $\G(\A_{f})$-modules:
	$$
		\tilde H^1(\G,\Q_{p})
		=
		\Hom_{\cts} (\Cong(\G),\Q_{p})_{\G(\A^{\gp})-\smooth}.
	$$
\end{thm}

\proof
Choose a level $K_f$ small enough so that
$\Gamma(K_f)$ is torsion-free.
By Corollary \ref{classgroup}, we have:
$$
	H^1(Y(K_f),\Z/p^s)
	=
	H^1_{\Group}(\Gamma(K_f),\Z/p^s).
$$
Elements of $H^1_{\Group}(\Gamma(K_f),\Z/p^s)$ are
group homomorphisms $\Gamma(K_f)\to \Z/p^s$.
Let $\tilde K_{f}$ be the preimage of $K_{f}$ in $\tilde \G(k)$;
this is equal to the profinite completion of $\Gamma(K_{f})$.
It follows that homomorphisms $\Gamma(K_{f})\to \Z/p^{s}$
correspond bijectively to continuous homomorphisms $\tilde K_{f}\to \Z/p^{s}$.
We therefore have:
$$
	H^1(Y(K_f),\Z/p^s)
	=
	H^1_{\disc}(\tilde K_f,\Z/p^s).
$$
We have an extension of profinite groups:
$$
	1
	\to
	\Cong(\G)
	\to
	\tilde{K}_f
	\to
	K_f
	\to
	1.
$$
This gives rise to a Hochschild--Serre spectral sequence
(Theorem \ref{hochschildserre}):
$$
	H^p_{\disc}(K_f,H^q_{\disc}(\Cong(\G),\Z/p^s))
	\implies
	H^{p+q}_{\disc}(\tilde K_f,\Z/p^s).
$$
From this we have an inflation-restriction sequence
containing the following terms:
\begin{equation}
	\label{infrest}
	0
	\to
	H^1_{\disc}(K_f,\Z/p^s)
	\to
	H^1(Y(K_f),\Z/p^s)
	\to
	H^1_{\disc}(\Cong(\G),\Z/p^s)^{K_f}
	\to
	H^2_{\disc}(K_f,\Z/p^s)
\end{equation}
The proof of the theorem consists of applying the functors $\limdKp$, $\limprojs$,
$-\otimes_{\Z_{p}}\Q_{p}$ and $\lim_{\to(K^{\gp})}$ to the sequence
(\ref{infrest}).
\medskip

\emph{Step 1.}
We first substitute $K_f=K_\gp K^\gp$, and
apply the functor $\limdKp$ to (\ref{infrest}).
We have by the K\"unneth formula:
$$
	\limdKp
	H^1_{\disc}(K_\gp K^{\gp},\Z/p^s)
	=
	H^{1}_{\disc}(K^{\gp},\Z/p^{s}),
$$
$$
	\limdKp
	H^{2}_{\disc}(K_\gp K^{\gp},\Z/p^s)
	=
	H^{2}_{\disc}(K^{\gp},\Z/p^{s}).
$$
By Lemma \ref{smooth} we have:
$$
	\limdKp
	H^1_{\disc}(\Cong(\G),\Z/p^s)^{K_{\gp}K^{\gp}}
	=
	H^1_{\disc}(\Cong(\G),\Z/p^s)^{K^{\gp}}.
$$
Since the functor $\limdKp$ is exact,
 the sequence remains exact:
\begin{equation}
	\label{infrest2}
	0
	\to
	H^{1}_{\disc}(K^{\gp},\Z/p^{s})
	\to
	\limdKp
	H^1(Y(K^{\gp}K_{\gp}),\Z/p^s)
	\to
	H^1_{\disc}(\Cong(\G),\Z/p^s)^{K^{\gp}}
	\to
	H^{2}_{\disc}(K^{\gp},\Z/p^{s}).
\end{equation}
\medskip

\emph{Interlude.}
Before going on, we make some restrictions on the
tame level $K^{\gp}$, and investigate the first and last
terms in the sequence (\ref{infrest2}).

We shall assume that the tame level $K^{\gp}$
 is a product of local factors:
$$
	K^{\gp}
	=
	\prod_{\gq\ne \gp}
	K_{\gq},
$$
where each $K_{\gq}$ is a compact open subgroup of $\G(k_{\gq})$.
Consider the following sets of finite primes of $k$:
$$
	S
	=
	\{\gq : \gq | p \hbox{ and } \gq \ne \gp \},
$$
$$
	T
	=
	\{
	\gq:
	\gq \not | p \hbox{ and } K_{\gq}\ne [K_{\gq},K_{\gq}]
	\}.
$$
Both these sets are finite.
We shall also assume from now on that for each prime $\gq\in T$,
the group $K_{\gq}$ is chosen small enough so that it is a
pro-$q$ group, where $q$ is the rational prime below $\gq$.
In particular, for each $\gq\in T$ we have for $n\ge 1$,
\begin{equation}
	\label{pf1}
	H^{n}_{\disc}(K_{\gq},\Z/p^{s})
	=
	0.
\end{equation}
We have a decomposition of $K^{\gp}$:
\begin{equation}
	\label{pf2}
	K^{\gp}
	=
	K_{S}
	\times
	K_{T}
	\times
	K^{S\cup T\cup\{\gp\}},
\end{equation}
where we are using the notation:
$$
	K_{S}
	=
	\prod_{\gq\in S} K_{\gq},\quad
	K^{S}
	=
	\prod_{\gq\not\in S} K_{\gq}.
$$
By the K\"unneth formula and (\ref{pf1}), (\ref{pf2}), we have:
\begin{equation}
	\label{pf3}
	H^{\bullet}_{\disc}(K^{\gp},\Z/p^{s})
	=
	H^{\bullet}_{\disc}(K_{S}K^{S\cup T\cup\{\gp\}},\Z/p^{s}).
\end{equation}
By assumption, the group $K^{S\cup T\cup\{\gp\}}$ is perfect,
so we have
\begin{equation}
	\label{pf4}
	H^{1}_{\disc}(K^{S\cup T\cup \{\gp\}},\Z/p^{s})
	=
	0.
\end{equation}
Again by the K\"unneth formula together with (\ref{pf3}), (\ref{pf4}),
 we have:
\begin{equation}
	\label{pf5}
	H^{1}_{\disc}(K^{\gp},\Z/p^{s})
	=
	H^{1}_{\disc}(K_{S},\Z/p^{s}).
\end{equation}
\begin{equation}
	\label{pf6}
	H^{2}_{\disc}(K^{\gp},\Z/p^{s})
	=
	H^{2}_{\disc}(K_{S},\Z/p^{s})
	\oplus
	H^{2}_{\disc}(K^{S\cup T\cup \{\gp\}},\Z/p^{s}).
\end{equation}
For each prime $\gq\notin S\cup T\cup\{\gp\}$, there is an open
normal pro-$q$ subgroup $L_{\gq}\subset K_{\gq}$.
We shall write $G(\gq)$ for the (finite) quotient group.
We therefore have a Hochschild--Serre spectral sequence:
$$
	H^{p}_{\Group}(G(\gq),H^{\gq}_{\disc}(L_{\gq},\Z/p^{s}))
	\implies
	H^{p+q}_{\disc}(K_{\gq},\Z/p^{s}).
$$
This spectral sequence degenerates: for $n\ge 1$ we have
$$
	H^{n}(L_{\gq},\Z/p^{s})
	=
	0.
$$
Hence,
\begin{equation}
	\label{pf7}
	H^{\bullet}_{\disc}(K_{\gq},\Z/p^{s})
	=
	H^{\bullet}_{\Group}(G(\gq),\Z/p^{s}),
	\qquad
	\gq\notin S\cup T\cup\{\gp\}.
\end{equation}
Since $G(\gq)$ is a finite perfect group, it has a universal central extension.
We shall write $\pi_{1}(G(\gq))$ for the kernel of this extension, i.e. the
Schur multiplier of $G(\gq)$.
By (\ref{pf7}) we have:
\begin{equation}
	\label{pf8}
	H^{2}_{\disc}(K_{\gq},\Z/p^{s})
	=
	\Hom_{\Group}(\pi_{1}(G(\gq)),\Z/p^{s}).
\end{equation}
By Corollary \ref{localglobal} and (\ref{pf8}) we have:
\begin{equation}
	\label{pf9}
	H^{2}_{\disc}(K^{S\cup T\cup \{\gp\}},\Z/p^{s})
	=
	\bigoplus_{\gq\notin S\cup T\cup \{\gp\}}
	\Hom_{\Group}(\pi_{1}(G(\gq)),\Z/p^{s}).
\end{equation}
From (\ref{pf6}) and (\ref{pf9}) we have:
\begin{equation}
	\label{pf9a}
	H^{2}_{\disc}(K^{\gp},\Z/p^{s})
	=
	H^{2}(K_{S},\Z/p^{s})
	\oplus
	\Hom_{\cts}(\pi_{1}^{S\cup T\cup \{\gp\}},\Z/p^{s}),
\end{equation}
where we are using the notation
$$
	\pi_{1}^{S\cup T\cup \{\gp\}}
	=
	\prod_{\gq\notin S\cup T\cup \{\gp\}}
	\pi_{1}(G(\gq)).
$$
The only property of $\pi_{1}^{S\cup T\cup \{\gp\}}$ which we require,
is that it is a product of finite groups, not depending on $s$.
\medskip

\emph{Step 2.}
We are now ready to apply the functor $\limprojs$
to the sequence (\ref{infrest2}).
To keep track of the exactness, we splice the sequence (\ref{infrest2}) into two:
\begin{equation}
	\label{infrest2a}
	0
	\to
	H^{1}_{\disc}(K_{S},\Z/p^{s})
	\to
	\limdKp
	H^1(Y(K^{\gp}K_{\gp}),\Z/p^s)
	\to
	A(s)
	\to
	0,
\end{equation}
\begin{equation}
	\label{infrest2b}
	0
	\to
	A(s)
	\to
	H^1_{\disc}(\Cong(\G),\Z/p^s)^{K^{\gp}}
	\to
	H^{2}_{\disc}(K^{\gp},\Z/p^{s}).
\end{equation}
\medskip

\emph{Step 2a.}
Applying the functor $\limprojs$ to (\ref{infrest2a}), we have
a long exact sequence:
\begin{equation}
	\label{pf10}
	0
	\to
	\limprojs
	H^{1}_{\disc}(K_{S},\Z/p^{s})
	\to
	\limprojs
	\limdKp
	H^1(Y(K^{\gp}K_{\gp}),\Z/p^s)
	\to
	\limprojs
	A(s)
	\to
	\limprojsd 
	H^{1}_{\disc}(K_{S},\Z/p^{s}).
\end{equation}
In order to calculate the individual terms in (\ref{pf10}),
we shall use the Eilenberg--Moore sequence
 (see Theorem \ref{eilenbergmoorethm}):
\begin{equation}
	\label{eilenbergmoore}
	0
	\to
	\limprojsd 
	H^{n-1}_{\disc}(K_{S},\Z/p^{s})
	\to
	H^{n}_{\cts}(K_{S},\Z_{p})
	\to
	\limprojs
	H^{n}_{\disc}(K_{S},\Z/p^{s})
	\to
	0.
\end{equation}
Taking $n=1$ in (\ref{eilenbergmoore}) we have
$$
	\limprojs
	H^{1}_{\disc}(K_{S},\Z/p^{s})
	=
	H^{1}_{\cts}(K_{S},\Z_{p}).
$$
Since $[K_{S},K_{S}]$ is open in $K_{S}$,
it follows that:
\begin{equation}
	\label{pf11}
	\limprojs
	H^{1}_{\disc}(K_{S},\Z/p^{s})
	=
	0.
\end{equation}
Also, since the groups $H^{1}_{\disc}(K_{S},\Z/p^{s})$
are all finite, it follows by Corollary \ref{finitegroups} that
\begin{equation}
	\label{pf12}
	\limprojsd 
	H^{1}_{\cts}(K_{S},\Z/p^{s})
	=
	0.
\end{equation}
Substituting (\ref{pf11}) and (\ref{pf12}) into (\ref{pf10}),
we get
\begin{equation}
	\label{pf13}
	\limprojs \limdKp
	H^1(Y(K^{\gp}K_{\gp}),\Z/p^s)
	=
	\limprojs
	A(s).
\end{equation}
\medskip

\emph{Step 2b.}
Applying the left-exact functor $\limprojs$ to (\ref{infrest2b})
and substituting (\ref{pf13}) we obtain the following exact sequence:
\begin{equation}
	\label{pf14}
	0
	\to
	\limprojs \limdKp
	H^1(Y(K^{\gp}K_{\gp}),\Z/p^s)
	\to
	\limprojs
	\left( H^1_{\disc}(\Cong(\G),\Z/p^s)^{K^{\gp}} \right)
	\to
	\limprojs
	H^{2}_{\disc}(K^{\gp},\Z/p^{s}).
\end{equation}
We shall investigate the second and third terms in this sequence further.

The functors $\limprojs$ and $-^{K^{\gp}}$ commute, so we have
\begin{equation}
	\label{pf15}
	\limprojs
	\left( H^1_{\disc}(\Cong(\G),\Z/p^s)^{K^{\gp}} \right)
	=
	\left(
	\limprojs
	H^1_{\disc}(\Cong(\G),\Z/p^{s})
	\right)^{K^{\gp}}.
\end{equation}
Again by the Eilenberg--Moore sequence (\ref{eilenbergmoore})
we have by (\ref{pf15}):
\begin{equation}
	\label{pf16}
	\limprojs
	\left( H^1_{\disc}(\Cong(\G),\Z/p^s)^{K^{\gp}} \right)
	=
	H^1_{\cts}(\Cong(\G),\Z_{p})^{K^{\gp}}.
\end{equation}
To calculate the third term in (\ref{pf14})
we shall use (\ref{pf9a}).
This shows that
\begin{equation}
	\label{pf17}
	\limprojs
	H^{2}_{\disc}(K^{\gp},\Z/p^{s})
	=
	\limprojs
	H^{2}_{\disc}(K_{S},\Z/p^{s})
	\oplus
	\limprojs
	\Hom_{\cts}(\pi_{1}^{S\cup T\cup\{\gp\}},\Z/p^{s}).
\end{equation}
Since $\pi_{1}^{S\cup T\cup\{\gp\}}$ is a product of
finite groups, it follows that
$$
	\limprojs
	\Hom_{\cts}(\pi_{1}^{S\cup T\cup\gp},\Z/p^{s})
	=
	0.
$$
Substituting this into (\ref{pf17}), we obtain:
\begin{equation}
	\label{pf18}
	\limprojs
	H^{2}_{\disc}(K^{\gp},\Z/p^{s})
	=
	\limprojs
	H^{2}_{\disc}(K_{S},\Z/p^{s}).
\end{equation}
Substituting (\ref{pf12}) into the Eilenberg--Moore sequence (\ref{eilenbergmoore}),
 we have:
\begin{equation}
	\label{pf19}
	\limprojs
	H_{\cts}^{2}(K_{S},\Z/p^{s})
	=
	H_{\cts}^{2}(K_{S},\Z_{p}).
\end{equation}
Substituting (\ref{pf19}) into (\ref{pf18}) we have:
$$
	\limprojs
	H^{2}_{\cts}(K^{\gp},\Z/p^{s})
	=
	H^{2}_{\cts}(K_{S},\Z_{p}).
$$
The sequence (\ref{pf14}) is therefore
\begin{equation}
	\label{infrest3}
	0
	\to
	\limprojs \limdKp  H^1(Y(K^{\gp}K_{\gp}),\Z/p^{s})
	\to
	H^{1}_{\cts}(\Cong(\G),\Z_{p})^{K^{\gp}}
	\to
	H^{2}_{\cts}(K_{S},\Z_{p}).
\end{equation}
\medskip

\emph{Step 3.}
We next apply the exact functor $-\otimes_{\Z_{p}}\Q_{p}$
to (\ref{infrest3}).
Note that since $K^{\gp}$ and $\Cong(\G)$ are compact,
we have
$$
	C^{\bullet}_{\cts}(K^{\gp},\Z_{p})\otimes_{\Z_{p}} \Q_{p}
	=
	C^{\bullet}_{\cts}(K^{\gp},\Q_{p}),
$$
$$
	C^{\bullet}_{\cts}(\Cong(\G),\Z_{p})\otimes_{\Z_{p}} \Q_{p}
	=
	C^{\bullet}_{\cts}(\Cong(\G),\Q_{p}).
$$
Furthermore, since $\Q_{p}$ is flat over $\Z_{p}$, we have
$$
	H^{\bullet}_{\cts}(K^{\gp},\Z_{p})\otimes_{\Z_{p}} \Q_{p}
	=
	H^{\bullet}_{\cts}(K^{\gp},\Q_{p}),
$$
$$
	H^{\bullet}_{\cts}(\Cong(\G),\Z_{p})\otimes_{\Z_{p}} \Q_{p}
	=
	H^{\bullet}_{\cts}(\Cong(\G),\Q_{p}).
$$
Since $H^{1}_{\cts}(\Cong(\G),\Z_{p})$ is torsion-free, we have
$$
	\left(H^{1}_{\cts}(\Cong(\G),\Z_{p})^{K^{\gp}}\right)\otimes_{\Z_{p}}\Q_{p}
	=
	H^{1}_{\cts}(\Cong(\G),\Q_{p})^{K^{\gp}}.
$$
Again, since $\Q_{p}$ is flat over $\Z_{p}$, we have an exact sequence:
\begin{equation}
	\label{infrest4}
	0
	\to
	\tilde H^1(K^{\gp},\Q_{p})
	\to
	H^{1}_{\cts}(\Cong(\G),\Q_{p})^{K^{\gp}}
	\to
	H^{2}_{\cts}(K_{S},\Q_{p}).
\end{equation}
\medskip

\emph{Step 4.}
Applying the exact functor $\lim_{\to\atop K^{\gp}}$ to (\ref{infrest4}),
 we have an exact sequence
\begin{equation}
	\label{infrest5}
	0
	\to
	\tilde H^1(\G,\Q_{p})
	\to
	H^{1}_{\cts}(\Cong(\G),\Q_{p})_{\G(\A^{\gp})-\smooth}
	\to
	H^{2}_{\st}(\G(k_{S}),\Q_{p}).
\end{equation}
As $\G(k_{S})$ is a $\Q_{p}$-analytic group,
 the stable cohomology may be expressed in terms of
 Lie algebra cohomology (using (\ref{stableLie})):
$$
	H^{2}_{\st}(\G(k_{S}),\Q_{p})
	=
	H^{2}_{\Lie}(\gg\otimes_{k}k_{S},\Q_{p}),
$$
where we are regarding $\gg\otimes_{k}k_{S}$ as a Lie algebra over $\Q_{p}$.
By Whitehead's second lemma (Theorem \ref{whitehead2})
we have
$$
	H^{2}_{\st}(\G(k_{S}),\Q_{p})
	=
	0.
$$
Hence
$$
	\tilde H^1(\G,\Q_{p})
	=
	H^{1}_{\cts}(\Cong(\G),\Q_{p})_{\G(\A^{\gp})-\smooth}.
$$
\sq
\medskip

\section{Some Examples}

\subsection{$\SL_{2}/\Q$}

Let $\G=\SL_{2}/\Q$.
Since $\gg$ is 3-dimensional, the spectral sequence has non-zero
terms only in columns $0$ to $3$.
Since arithmetic subgroups have virtual cohomological dimension $1$,
it follows that $\tilde H^{n}=0$ for $n>1$.
Taking $W$ to be the trivial representation,
the $E_{2}$ sheet of the spectral sequence is as follows:
$$
	E_{2}^{\bullet,\bullet}\quad : \quad
	\begin{array}{cccc}
		\tilde H^{1}(\G,E)^{\gg}_{\la}& E & 0 & 0\medskip\\
		E & 0 & 0 & E
	\end{array}
$$
The connection map $E_{2}^{1,1}\to E_{2}^{3,0}$ is an isomorphism,
and the spectral sequence stabilizes at $E_{3}$ as follows:
$$
	E_{3}^{\bullet,\bullet}\quad : \quad
	\begin{array}{cccc}
		\tilde H^{1}(\G,E)^{\gg}_{\la}& 0 & 0 & 0\medskip\\
		E & 0 & 0 & 0
	\end{array}
$$

\subsection{$\SL_{1}(D)$ for an indefinite quaternion algebra $D$}

Let $k$ be a totally real field and let $D$ be a quaternion algebra over $k$,
which is indefinite at exactly one real place of $k$.
We shall consider the group $\G(-)=\SL_{1}(D\otimes_{k}-)$ over $k$.
Arithmetic subgroups of $\G$ have virtual cohomological dimension $2$,
so we have classical cohomology groups in dimensions $0$, $1$ and $2$.
In dimensions $0$ and $2$ these are given by constants,
 and are $1$-dimensional.
On the other hand it is easy to show that $\tilde H^{2}(\G,\Q_{p})=0$.
The $E_{2}$ sheet of the spectral sequence is as follows:
$$
	E_{2}^{\bullet,\bullet}
	\qquad:
	\qquad
	\begin{array}{cccc}
		\tilde H^{1}(\G,E)^{\gg}_{\la}& E^{2} & 0 & 0\medskip\\
		E & 0 & 0 & E
	\end{array}
$$
The connection map $E_{2}^{1,1}\to E_{2}^{3,0}$ is surjective,
and the spectral sequence stabilizes at $E_{3}$ as follows:
$$
	E_{3}^{\bullet,\bullet}
	\qquad:
	\qquad
	\begin{array}{cccc}
		\tilde H^{1}(\G,E)^{\gg}_{\la}& E & 0 & 0\medskip\\
		E & 0 & 0 & 0
	\end{array}
$$

\subsection{$\SL_{2}/k$ for $k$ real quadratic}

Let $k$ be a real quadratic field and consider the group $\G=\SL_{2}/k$.
The non-zero classical cohomology groups are the following:
\begin{eqnarray*}
	H^{0}_{\cl}(\G,W)
	&=&
	W^{\G},\\
	H^{2}_{\cl}(\G,W)
	&&
	\hbox{infinite dimensional}.
\end{eqnarray*}
It is known in this case (see \cite{serre3})
 that the congruence kernel of $\G$ is trivial.
We therefore have $\tilde H^{1}(\G,E)=0$,
and we can also show that $\tilde H^{3}(\G,E)=0$.
Therefore the weak Emerton criterion holds in dimension 2.
We also have $H^{3}(\gg,K_{\infty},\C)=0$.
Therefore we may apply Theorem \ref{variation2}
to the eigenvariety $E(2,K^{\gp})$.
The $E_{2}$-sheet of the spectral sequence is as follows:
$$
	E_{2}^{\bullet,\bullet}
	\qquad:
	\qquad
	\begin{array}{cccc}
		\tilde H^{2}(\G,E)^{\gg}_{\la}& 0 & 0 & 0\medskip\\
		0& 0 & 0 & 0\medskip\\
		E & 0 & 0 & E
	\end{array}
$$
The map $\tilde H^{2}(\G,E)_{\la}^{\gg}\to E$ in the $E_{3}$-sheet
 is surjective, and the spectral sequence stabilizes at the $E_{4}$-sheet:
$$
	E_{4}^{\bullet,\bullet}
	\qquad:
	\qquad
	\begin{array}{cccc}
		H^{2}_{\cl}(\G,E)& 0 & 0 & 0\medskip\\
		0& 0 & 0 & 0\medskip\\
		E & 0 & 0 & 0
	\end{array}
$$

\subsection{$\SL_{3}/\Q$}

Arithmetic subgroups of $\SL_{3}(\Q)$ have virtual cohomological dimension $4$, as the symmetric space is $5$-dimensional.
We have the following non-zero classical cohomology groups:
\begin{eqnarray*}
	H^{0}_{\cl}(\G,W)
	&=&
	W^{\G},\\
	H^{2}_{\cl}(\G,W)
	&=&
	\hbox{infinite dimensional},\\
	H^{3}_{\cl}(\G,W)
	&=&
	\hbox{infinite dimensional}.
\end{eqnarray*}
It was shown in \cite{basslazardserre} that the congruence kernel is trivial.
Hence the weak Emerton criterion holds in dimension 2,
and in fact the only non-zero Banach space representations are:
\begin{eqnarray*}
	\tilde H^{0}(\G,E)
	&=&
	E,\\
	\tilde H^{2}(\G,E)
	&=&
	\hbox{infinite dimensional},\\
	\tilde H^{3}(\G,E)
	&=&
	\hbox{infinite dimensional}.
\end{eqnarray*}
Furthermore, $H^{3}_{\relLie}(\gg,K_{\infty},\C)=0$.
We may therefore apply Theorem \ref{variation2}
to the eigenvariety $E(2,K^{\gp})$.
One can use Poincar\'e duality to construct an eigenvariety
 interpolating $H^{3}_{\cl}$.

The author has not been able to calculate
all of the terms of the spectral sequence.
However the $E_{2}$-sheet is as follows:
$$
	E_{2}^{\bullet,\bullet}
	\qquad:
	\qquad
	\begin{array}{ccccccccc}
		\tilde H^{3}(\G,E)^{\gg}_{\la}& ? & ? & ?& ? & ? & ? & 0 & 0\medskip\\
		\tilde H^{2}(\G,E)^{\gg}_{\la}& \Ext^{1}_{\gg}(E,\tilde H^{2}(\G,E)_{\la}) & ? & ?& ? & ? & ? & ? & ?\medskip\\
		0& 0 & 0 & 0 & 0 & 0 & 0 & 0 & 0\medskip\\
		E & 0 & 0 & E & 0 & E & 0 & 0 & E
	\end{array}
$$
This is stable by the $E_{5}$-sheet, and most things are known:
$$
	E_{5}^{\bullet,\bullet}
	\qquad:
	\qquad
	\begin{array}{ccccccccc}
		?& 0 & 0 & 0& 0 & 0 & 0 & 0 & 0\medskip\\
		H^{2}_{\cl}(\G,E)& ? & 0 & 0& 0 & 0 & 0 & 0 & 0\medskip\\
		0& 0 & 0 & 0 & 0 & 0 & 0 & 0 & 0\medskip\\
		E & 0 & 0 & 0 & 0 & 0 & 0 & 0 & 0
	\end{array}
$$

\subsection{$\Sp_{4}/\Q$}

Arithmetic subgroups of $\Sp_{4}(\Q)$ have cohomological dimension $5$, as the symmetric space is $6$-dimensional.
It was shown in \cite{bassmilnorserre}
 that the congruence kernel is trivial.
Furthermore $H^{3}(\gg, K_{\infty},\C)=0$.
We may therefore apply Theorem \ref{variation2}
to give a construction of the $H^{2}$-eigencurve.
By Poincar\'e duality, it is also possible to construct a reasonable
 $H^{4}$-eigenvariety.

\subsection{$\Spin(2,l)$ ($l\ge 3$)}

Let $L$ be a $\Z$-lattice equipped with a quadratic form
of signature $(2,l)$ with $l\ge 3$.
We let $\G/\Q$ be the corresponding Spin group.
This has real rank $2$, and the corresponding symmetric space has dimension $2l$.
The congruence kernel was shown to be trivial
for such groups by Kneser \cite{kneserspin}.
Hence $\G$ satisfies the weak Emerton criterion in dimension $2$.
It turns out that $H^{3}(\gg,K_{\infty},\C)=0$,
 so we may apply Theorem \ref{variation2} to $E(2,K^{\gp})$.


\begin{thebibliography}{10}
\bibitem{ashpollackstevens}
A.~Ash, D.~Pollack, and G.~Stevens.
\newblock Rigidity of $p$-adic cohomology classes of congruence subgroups of
  $\GL(n,\Z)$.
\newblock {\em preprint}, 2007.

\bibitem{ashstevens}
A.~Ash and G.~Stevens.
\newblock {$p$}-adic deformations of cohomology classes of subgroups of {${\rm
  GL}(n,{\bf Z})$}.
\newblock {\em Collect. Math.}, 48 (1--2) 1--30, 1997.
\newblock Journ\'ees Arithm\'etiques (Barcelona, 1995).

\bibitem{ashstevens2}
A.~Ash and G.~Stevens.
\newblock $p$-adic deformations of cohomology classes of subgroups of $\GL(N,\Z)$: the non-ordinary case.
\newblock {\em preprint}, April 2000.

\bibitem{basslazardserre}
H.~Bass, M.~Lazard, and J.-P. Serre.
\newblock Sous-groupes d'indice fini dans {${\bf SL}(n,\,{\bf Z})$}.
\newblock {\em Bull. Amer. Math. Soc.}, 70:385--392, 1964.

\bibitem{bassmilnorserre}
H.~Bass, J.~Milnor, and J.-P. Serre.
\newblock Solution of the congruence subgroup problem for $\SL_n$  ($n\ge 3$) and
  $\Sp_{2n}$ ($n \ge 2$).
\newblock {\em Publications math{\'e}matiques de l'I.H.{\'E}.S.}, 33:59--137,
  1967.

\bibitem{blasiusfrankegrunewald}
D.~Blasius, J.~Franke, and F.~Grunewald.
\newblock Cohomology of {$S$}-arithmetic subgroups in the number field case.
\newblock {\em Invent. Math.}, 116(1-3):75--93, 1994.

\bibitem{borel2}
A.~Borel.
\newblock Cohomologie de sous-groupes discrets et repr\'esentations de groupes
  semi-simples.
\newblock In {\em Colloque ``Analyse et Topologie'' en l'Honneur de Henri
  Cartan (Orsay, 1974)}, pages 73--112. Ast\'erisque, No. 32--33. Soc. Math.
  France, Paris, 1976.

\bibitem{borelbook}
A.~Borel.
\newblock {\em Linear algebraic groups}, volume 126 of {\em Graduate Texts in
  Mathematics}.
\newblock Springer-Verlag, New York, second edition, 1991.

\bibitem{borelwallach}
A.~Borel and N.~R. Wallach.
\newblock {\em Continuous cohomology, discrete subgroups, and representations
  of reductive groups}, vol. 67 of {\em Mathematical surveys and
  monographs}.
\newblock American Mathematical Society, Providence, R.I., 2nd edition,
  2000.

\bibitem{buzzard2004}
K.~Buzzard.
\newblock On {$p$}-adic families of automorphic forms.
\newblock In {\em Modular curves and abelian varieties}, volume 224 of {\em
  Progr. Math.}, pages 23--44. Birkh\"auser, Basel, 2004.


\bibitem{buzzard3}
K.~Buzzard.
\newblock Eigenvarieties.
\newblock {\em preprint}, 2006.



\bibitem{cartier}
P.~Cartier.
\newblock Representations of {$p$}-adic groups: a survey.
\newblock In {\em Automorphic forms, representations and $L$-functions (Proc.
  Sympos. Pure Math., Oregon State Univ., Corvallis, Ore., 1977), Part 1},
  Proc. Sympos. Pure Math., XXXIII, pages 111--155. Amer. Math. Soc.,
  Providence, R.I., 1979.


\bibitem{chenevier}
G.~Chenevier.
\newblock Familles {$p$}-adiques de formes automorphes pour {${\rm GL}\sb n$}.
\newblock {\em J. Reine Angew. Math.}, 570:143--217, 2004.

\bibitem{coleman}
R.~F. Coleman.
\newblock {$p$}-adic {B}anach spaces and families of modular forms.
\newblock {\em Invent. Math.}, 127(3):417--479, 1997.

\bibitem{colemanmazur}
R.~Coleman and B.~Mazur.
\newblock The eigencurve.
\newblock In {\em Galois representations in arithmetic algebraic geometry
  (Durham, 1996)}, volume 254 of {\em London Math. Soc. Lecture Note Ser.},
  pages 1--113. Cambridge Univ. Press, Cambridge, 1998.

\bibitem{emertonAnalytic}
M.~Emerton.
\newblock {\em Locally analytic vectors in representations of locally $p$-adic
  analytic groups}.
\newblock preprint, 2004.

\bibitem{emertonJacq}
M.~Emerton.
\newblock Jacquet modules of locally analytic representations of {$p$}-adic
  reductive groups. {I}. {C}onstruction and first properties.
\newblock {\em Ann. Sci. \'Ecole Norm. Sup. (4)}, 39(5):775--839, 2006.

\bibitem{emerton}
M.~Emerton.
\newblock On the interpolation of systems of eigenvalues attached to
  automorphic {H}ecke eigenforms.
\newblock {\em Invent. Math.}, 164(1):1--84, 2006.

\bibitem{franke}
J.~Franke.
\newblock Harmonic analysis in weighted {$L\sb 2$}-spaces.
\newblock {\em Ann. Sci. \'Ecole Norm. Sup. (4)}, 31(2):181--279, 1998.

\bibitem{hida}
H.~Hida.
\newblock {\em Elementary theory of {$L$}-functions and {E}isenstein series},
  volume~26 of {\em London Mathematical Society Student Texts}.
\newblock Cambridge University Press, Cambridge, 1993.

\bibitem{kassaei-2004}
P.~L. Kassaei.
\newblock {$\mathscr P$}-adic modular forms over {S}himura curves over totally real
  fields.
\newblock {\em Compos. Math.}, 140(2):359--395, 2004.


\bibitem{MR2154369}
M.~Kisin and K.~F. Lai.
\newblock Overconvergent {H}ilbert modular forms.
\newblock {\em Amer. J. Math.}, 127(4):735--783, 2005.

\bibitem{kneser1}
M.~Kneser.
\newblock Starke {A}pproximation in algebraischen {G}ruppen. {I}.
\newblock {\em J. Reine Angew. Math.}, 218:190--203, 1965.

\bibitem{kneser-hasseprinciple}
M.~Kneser.
\newblock Hasse principle for {$H\sp{1}$} of simply connected groups.
\newblock In {\em Algebraic Groups and Discontinuous Subgroups (Proc. Sympos.
  Pure Math., Boulder, Colo., 1965)}, pages 159--163. Amer. Math. Soc.,
  Providence, R.I., 1966.

\bibitem{kneser}
M.~Kneser.
\newblock Strong approximation.
\newblock In {\em Algebraic Groups and Discontinuous Subgroups (Proc. Sympos.
  Pure Math., Boulder, Colo., 1965)}, pages 187--196. Amer. Math. Soc.,
  Providence, R.I., 1966.

\bibitem{kneserspin}
M.~Kneser.
\newblock Normalteiler ganzzahliger {S}pingruppen.
\newblock {\em J. Reine Angew. Math.}, 311/312:191--214, 1979.

\bibitem{neukirch}
J.~Neukirch, A.~Schmidt, and K.~Wingberg.
\newblock {\em Cohomology of number fields}, volume 323 of {\em Grundlehren der
  Mathematischen Wissenschaften [Fundamental Principles of Mathematical
  Sciences]}.
\newblock Springer-Verlag, Berlin, 2000.


\bibitem{platonov}
V.~P. Platonov.
\newblock The problem of strong approximation and the {K}neser-{T}its
  hypothesis for algebraic groups.
\newblock {\em Izv. Akad. Nauk SSSR Ser. Mat.}, 33:1211--1219, 1969.

\bibitem{prasad}
G.~Prasad.
\newblock Strong approximation for semi-simple groups over function fields.
\newblock {\em Ann. of Math. (2)}, 105(3):553--572, 1977.

\bibitem{raghunathan}
M.~S. Raghunathan.
\newblock The congruence subgroup problem.
\newblock {\em Proc. Indian Acad. Sci. Math. Sci.}, 114(4):299--308, 2004.

\bibitem{rapinchuk}
A.~S. Rapinchuk.
\newblock Congruence subgroup problem for algebraic groups: old and new.
\newblock {\em Ast\'erisque}, (209):11, 73--84, 1992.
\newblock Journ\'ees Arithm\'etiques, 1991 (Geneva).

\bibitem{schneider}
P.~Schneider.
\newblock Basic notions of rigid analytic geometry.
\newblock In {\em Galois representations in arithmetic algebraic geometry
  (Durham, 1996)}, volume 254 of {\em London Math. Soc. Lecture Note Ser.},
  pages 369--378. Cambridge Univ. Press, Cambridge, 1998.

\bibitem{schneider-teitelbaum}
P.~Schneider and J.~Teitelbaum.
\newblock Banach space representations and {I}wasawa theory.
\newblock {\em Israel J. Math.}, 127:359--380, 2002.

\bibitem{schneider-teitelbaum2}
P.~Schneider and J.~Teitelbaum.
\newblock Locally analytic distributions and {$p$}-adic representation theory,
  with applications to {${\rm GL}\sb 2$}.
\newblock {\em J. Amer. Math. Soc.}, 15(2):443--468 (electronic), 2002.

\bibitem{serre}
J.-P. Serre.
\newblock {\em Cohomologie galoisienne}, volume~5 of {\em Lecture Notes in
  Mathematics}.
\newblock Springer-Verlag, Berlin, fifth edition, 1994.

\bibitem{serre2}
J.-P. Serre.
\newblock Cohomologie des groupes discrets.
\newblock In {\em Prospects in mathematics (Proc. Sympos., Princeton Univ.,
  Princeton, N.J., 1970)}, pages 77--169. Ann. of Math. Studies, No. 70.
  Princeton Univ. Press, Princeton, N.J., 1971.

\bibitem{serre3}
J.-P. Serre.
\newblock Le probl\`eme des groupes de congruence pour $\mathrm{SL}_2$.
\newblock {\em Ann. of Math. (2)}, 92:489--527, 1970.

\bibitem{tits}
J.~Tits.
\newblock Reductive groups over local fields.
\newblock In {\em Automorphic forms, representations and $L$-functions (Proc.
  Sympos. Pure Math., Oregon State Univ., Corvallis, Ore., 1977), Part 1},
  Proc. Sympos. Pure Math., XXXIII, pages 29--69. Amer. Math. Soc., Providence,
  R.I., 1979.

\bibitem{weibel}
C.~A. Weibel.
\newblock {\em An introduction to homological algebra}, volume~38 of {\em
  Cambridge Studies in Advanced Mathematics}.
\newblock Cambridge University Press, Cambridge, 1994.

\end{thebibliography}
\end{document}